\documentclass[journal]{IEEEtran}
\usepackage{epsfig}
\usepackage{graphicx}
\usepackage{amsthm,epsfig}
\usepackage{tikz}
\usepackage{amssymb}
\usepackage{amsmath}
\usepackage{mathrsfs}
\usetikzlibrary{shapes,decorations}
\tikzstyle{every picture}+=[remember picture]

\newcommand{\cons}[1][]{\ensuremath{\alpha\ifthenelse{\equal{#1}{}}{}{^{(#1)}}}\xspace}

\newcommand{\F}{{\mathcal{F}}}

\newcommand{\e}{{\bf e}}

\newcommand{\mS}{{\mathcal{S}}}
\newcommand{\ra}{\rightarrow}

\newtheorem{theorem}{Theorem}

\newtheorem{prop}{Proposition}

\begin{document}

\title{A Mean-field Approach for Controlling Singularly Perturbed Multi-population SIS Epidemics}

  \author{P. T. Akhil$^1$, Eitan Altman$^2$, Rajesh Sundaresan$^{1,3}$\thanks{This work was supported in part by the Defence Research and Development Organisation [grant no. DRDO0667] under the DRDO-IISc Frontiers Research Programme and in part by the Indo-French Centre for Applied Mathematics.}
      \\$^1$ECE Department, IISc Bangalore, India.
      \\$^2$INRIA Sophia-Antipolis NEO team; University of Cote d'Azur, France; LIA, University of Avignon, France; and LINCS, Paris.
      \\$^3$Robert Bosch Centre for Cyber-Physical Systems, IISc Bangalore, India.
      }
      \date{}

\maketitle

\begin{abstract}
We consider a multi-population epidemic model with one or more (almost) isolated communities and one mobile community. Each of the isolated communities has contact within itself and, in addition, contact with the outside world but only through the mobile community. The contact rate between the mobile community and the other communities is assumed to be controlled. We first derive a multidimensional ordinary differential equation (ODE) as a mean-field fluid approximation to the process of the number of infected nodes, after appropriate scaling. We show that the approximation becomes tight as the sizes of the communities grow. We then use a singular perturbation approach to reduce the dimension of the ODE and identify an optimal control policy on this system over a fixed time horizon via Pontryagin's minimum principle. We then show that this policy is close to optimal, within a certain class, on the original problem for large enough communities. From a phenomenological perspective, we show that the epidemic may sustain in time in all communities (and thus the system has a nontrivial metastable regime) even though in the  absence of the mobile nodes the epidemic would die out quickly within each of the isolated communities.
\end{abstract}

\section{Introduction}\label{sec:intro}
In recent years there has been a growing interest in applying epidemics-related control techniques to networks. The control may have as objective either to limit the propagation of content when the content consists of malware or e-viruses or, on the contrary, to help propagate content when it consists of advertisements, news, entertainment, sports events, etc. In this paper, we consider an information diffusion problem wherein the objective is to obtain a good tradeoff between the information spread in the network and in the use of system resources. The spread of information closely resembles an epidemic spread. Epidemiological models are therefore used in modeling this information spread. Our work uses the susceptible-infected-susceptible (SIS) epidemic model. A member that has a copy of the information or content is said to be {\em infected} and a member that does not have the copy of the content is said to be {\em susceptible}. When two members come in contact, one infected and the other susceptible, the former transmits a copy of the content to the latter, and the latter gets infected. An infected member may spontaneously get rid of the content, a phenomenon that we call {\em curing}, to become susceptible again. Our social network model has a fraction of ``influential'' members whose interaction with other members can be controlled. The spread of information is achieved by giving incentives to these members to actively spread it, but this imposes a cost to the campaigner. We assume that initially only a certain fraction of members possess the information. We aim to maximize the information spread at the end of a finite time horizon subject to the costs incurred in controlling the influential members' interactions.
Some interesting features of our social network model are as follows.
\begin{itemize}
\item Our social network model departs from other models in that it consists of a \emph{mobile} community of influential members and a finite number of \emph{isolated} communities. The members of each isolated community interact among themselves and also with the members of the mobile community. But there is no {\em direct} interaction between the members of different isolated communities. Also, the members of the mobile community interact among themselves.

    Interactions in social media often happen in almost closed communities with a few influential members interacting across groups. The mobile community models these influential members. The same model is also applicable in epidemic settings. The mobile community can be seen as tourists traveling across the globe and individual countries as isolated communities. In this context, however, the campaigner may be a health-care worker and his goal may be to contain an epidemic.

\item The interactions between members/nodes within an isolated community happen at a faster timescale compared to the interactions between members/nodes within the mobile community. Consequently, the isolated nodes are quick to get infected and quick to recover compared to the mobile nodes. In an opinion dynamics setup, the mobile nodes model stubborn individuals who are difficult to influence and isolated nodes model those that are more easily influenced or are members of a much more interactive network.

\item The campaigner has control only over the rates of interactions of the mobile community. This is the situation when the campaigner is able to give incentives only to members of the mobile community. Larger the rates of interactions of the mobile community, greater the cost to the campaigner. In epidemic settings, the campaigner incurs cost while restraining the interactions of the mobile community with the isolated communities.
\end{itemize}

In this paper, we seek to answer the following two questions.
\begin{itemize}
 \item Can the cross-community interactions help sustain infection in all the communities for cases wherein the infection would have died out otherwise? This would then be a nontrivial metastable regime of sustained infection.
 \item What is the resource allocation strategy that maximizes the information spread?
\end{itemize}
In order to gain understanding of the evolution of the infection, we derive the fluid limit of the system as the population of each community is scaled to infinity. The fluid limit is a two timescale dynamical system: the dynamics of infection in the isolated communities happens at a faster timescale compared to that of the mobile community. We further drive the timescale separation to be sufficiently large to reduce the fluid limit to a uni-dimensional dynamical system wherein the infection dynamics of the mobile community sees only the equilibrated values of the fast timescale infection dynamics corresponding to the isolated communities. We show that the evolution of the empirical distribution of the infected population in the original model and in the reduced dynamical system are close to each other. Hence we use the reduced dynamical system to study the survival of infection in our social network model.

\begin{figure}
  \centering
  \includegraphics[scale=.4]{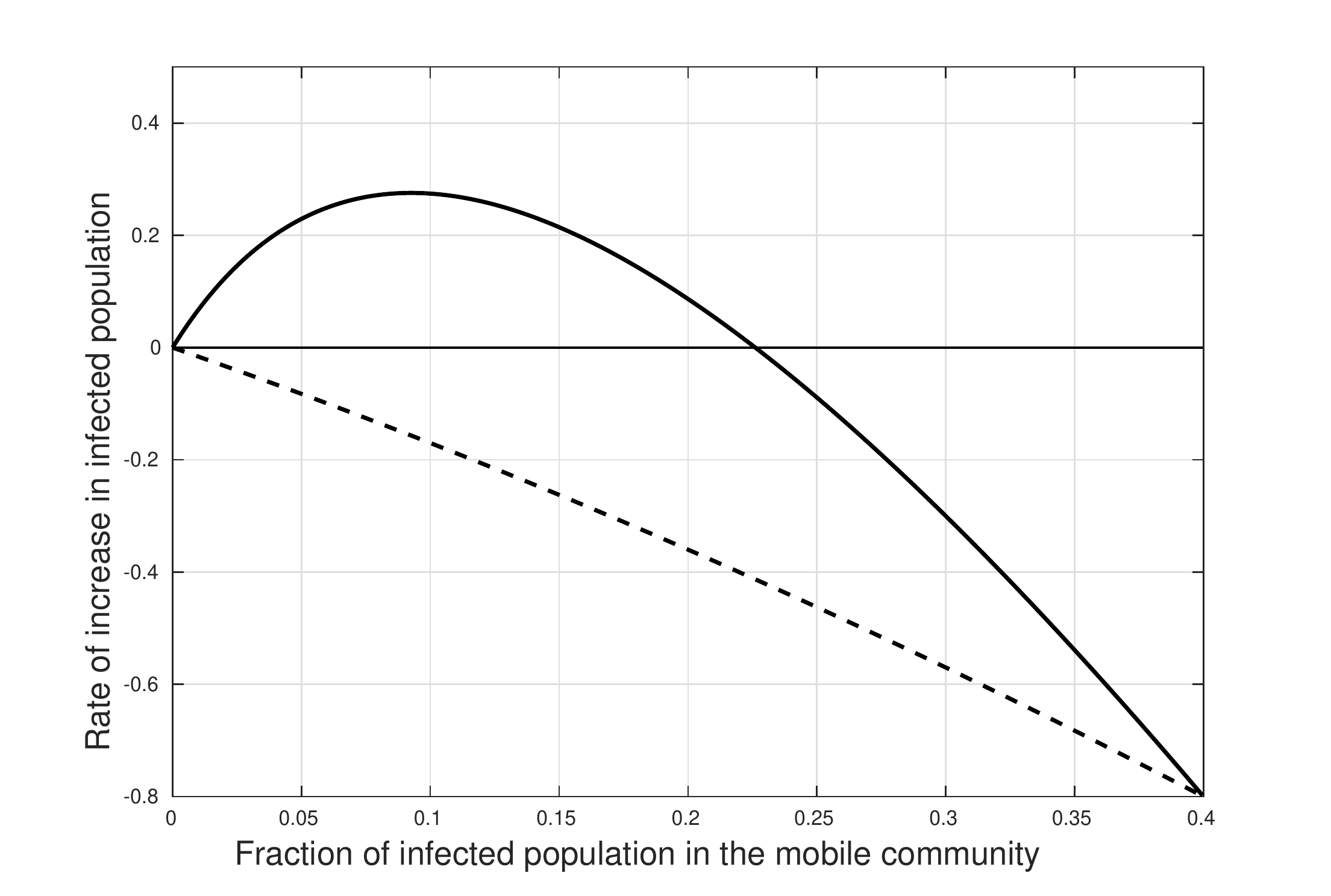}
  \caption{Cross-community interactions can help sustain infection. The picture is for the case of a population model with two isolated communities and a mobile community. The solid line corresponds to the case when cross-community interactions are present whereas the dashed line corresponds to the case of zero cross-community interactions. The zero-crossing points of the graphs show the equilibrium infection level.}
  \label{fig:network-effect-helps}
\end{figure}

We first consider the case when the rate of infection is less than the rate of recovery for each community considered in isolation. When the cross-community interactions are absent, the infection dies out, a fact that can be verified analytically and can also be gleaned from Figure \ref{fig:network-effect-helps}. Let us use the phrase ``infection level'' to refer to the fraction of infected members in a community. Figure \ref{fig:network-effect-helps} shows the rate of change of infection level in the mobile community as a function of the infection level in that community, in the reduced dynamical system. The dashed line plots this for the case when the cross-community interactions are absent. In this case, the rate of change is negative whenever a strictly positive fraction of the mobile community is infected. Thus the equilibrium is at $0$ and the infection dies out. The solid line represents the case when cross-community interactions are present. It is easily seen that the equilibrium value of the infection level in the mobile community corresponds to the zero crossing point of the graph and one of them is strictly positive. We will soon see that the infection sustains in this case.
However it must be noted that there exist interaction rates for which the infection dies out even in the presence of cross-community interactions. We characterize the set of interaction rates for which the infection sustains. There is more discussion on this phenomenon later in this section.

We next consider the optimal resource allocation strategy to maximize information spread. We define a cost function that has two components 1) a running cost that depends on both the instantaneous control applied and the instantaneous fraction of infected members, and 2) a terminal cost that depends on the fraction of infected members at the end of the finite time horizon. We study optimal control in the following steps.

\begin{enumerate}
\item \label{item:fluid-limit} We compare the costs on the original finite population system and the fluid limit. We show that for an identical cost function and identical path-wise control, the cost incurred on the original system and on the fluid limit are close to each other.
\item We then show that, for a specific cost function with a running cost linearly increasing in control (the intensity of interaction between the mobile community and the isolated communities) and a terminal cost linearly decreasing in the fraction of infected nodes of the mobile community, and for the same path-wise control with a fixed maximum number of discontinuities, the costs incurred on the fluid model and on the reduced dynamical system are close to each other.
\item Using the above results, we then show  that the cost of control on the original system and on the reduced dynamical system are close to each other.
\item Finally, we apply Pontryagin's minimum principle to show that a bang-bang control is optimal for the reduced fluid limit dynamical system. Hence that same control is nearly optimal on the finite population system, within the class of all policies that yield control-process sample paths with a fixed maximum number of discontinuities.
\end{enumerate}

Our choice of cost function for the above analysis corresponds to the case of maximization of information spread. In the case of disease control, we would minimize epidemic spread, and our analysis for maximizing information spread can be easily adapted to minimize epidemic spread.

From a technical standpoint, the difficulties to surmount are the following. First, one may view the finite but large population (say of size $n$) system as a small noise perturbation of a deterministic fluid limit. As we will see, the timescale separation parameter $\varepsilon$ will be taken to be $\varepsilon = C / (\log n)$ for some $0 < C < \infty$. As $n \rightarrow \infty$ one then encounters simultaneous small noise phenomena in both slow and fast timescales with the additional complication that the timescale separation $\varepsilon$ is controlled by the same parameter\footnote{A similar situation, simultaneous small noise phenomena in slow and fast timescales time scales with timescale separation and both small noises parameterized by the same $\varepsilon$, is studied in a system of diffusions without control \cite{athreya2018simultaneous}. There too a similar inverse logarithmic relation between the timescale separation $\varepsilon$ and the small noise level $1/n$ is required to ensure that the fast timescale variables converge to certain desired points.}. There is however sufficient regularity on the fast timescale subsystem that it concentrates near an equilibrium which can be explicitly identified. This is then exploited to understand the limiting system. Second, the limiting reduced dynamical system is one for which a direct argument does not immediately provide existence of an optimal control policy. Instead, we prove existence indirectly by comparison with a simpler system for which existence of an optimal control can be established, and then by arguing that the optimal control on the simpler system applies on the original system.

From a phenomenological standpoint, what is interesting is the following. There is one eventual absorbing state for the system, which is the state when the epidemic has completely died out. However, with interaction between communities and time scale separation, the entry into the absorbing state can be delayed well beyond any fixed and finite horizon, and the system can be maintained in a metastable equilibrium of sustained infection over any finite horizon. This phenomenon has been studied before. See \cite{ottaviano2017optimal} and references therein. For the model in \cite{ottaviano2017optimal} with infection rate to curing rate ratios below the so-called ``epidemic threshold'', the epidemic dies out quickly; if above the threshold, the infection sustains for a long time. There is quite a bit of focus in the literature on determining the epidemic threshold, again see \cite{ottaviano2017optimal} and references therein. We provide a sufficient condition for infection to sustain in Section \ref{sec:infection-sustains}. When there is a cost to maintain the interaction, a bang-bang control is optimal for a specific choice of cost functions. The specific cost function given later in (\ref{eqn:cost-function}) is mainly for illustration and extensions to more general costs should be possible and is left for the future.

{\it Related works:} The design of optimal control of epidemic spread has been studied in the context of disease/infection control in human networks \cite{morton1974optimal}, \cite{behncke2000optimal}, information or opinion dynamics in social networks \cite{karnik2012optimal}, \cite{kandhway2014run}, \cite{altman2010optimal}, security in mobile networks \cite{khouzani2012optimal}, etc. Kumar et al. \cite{kumar2018influenzing} studied the optimal timing of external influence on opinion dynamics in a homogeneous voter model. Karnik et al. \cite{karnik2012optimal} and Kandhway and Kuri \cite{kandhway2014run} designed optimal control for maximizing information spread in the following cases: a) an SIR epidemic model where the intensity of recruiting spreaders is controlled and b) an SIS epidemic model where, in addition to the recruitment of direct spreaders, a word-of-mouth spreading is also controlled. The works \cite{karnik2012optimal} and \cite{kandhway2014run} used cost functions that are linear and quadratic in the control variable, respectively. Khouzani et al. \cite{khouzani2012optimal} designed optimal control for the spread of security patches in mobile wireless networks. The fraction of nodes that actively spread security patches and the patch transmission rate are controlled. Altman et al. \cite{altman2010optimal} found the optimal strategy for activation and transmission power control of nodes in delay tolerant networks. These works used a general cost function compared to the linear and quadratic functions in \cite{karnik2012optimal} and \cite{kandhway2014run}. The aforementioned works considered the underlying network to be homogeneous where nodes mix uniformly.

Most practical scenarios, however, have an underlying network that is not homogeneous. Colizza et al. \cite{colizza2006role} modeled the spread of a global epidemic on a network with nodes denoting airports of major cities and edge weights accounting for the passenger flow between them. The resulting network is highly heterogeneous in traffic pattern between cities. Becker and Dietz \cite{becker1995effect} and Ball et al. \cite{ball2004stochastic} modeled the household structure of human population with different rates for within-household interactions and between-household interactions for a study to determine the fraction of population to be vaccinated to prevent an epidemic. Aditya et al. \cite{gopalan2011random} studied epidemic spreading in a large static network where the epidemic spread by means of a set of virtual mobile agents that can infect any node in the network in addition to the spread via usual node to node interactions.
They show that a small number of virtual mobile agents speed up the spread of epidemic to all nodes. Pellegrini et al. \cite{de2010optimal} studied optimal control of information spread in delay tolerant networks having multiple classes of nodes. The control variable is the probability of forwarding a message between any two classes of nodes. However, \cite{gopalan2011random} and \cite{de2010optimal} used a monotone SI model of epidemic spread unlike our SIS model. Kandhway and Kuri \cite{kandhway2016campaigning} extended the work in \cite{kandhway2014run} to a general graph model, but the optimal control was numerically obtained. Ottaviano et al. \cite{ottaviano2017optimal} formulated the optimal node-specific control policy for epidemic control in a heterogeneous network as a solution to an SDP. When the network has a community structure like the model considered in this paper, they achieve dimensionality reduction by modeling the network as a graph with equitable partition. Our problem is a finite time horizon epidemic control seeking a good trade-off between the cost of containing the epidemic and number of infected nodes after a fixed finite duration unlike \cite{ottaviano2017optimal} where the objective is to minimize curing costs subject to the condition that the epidemic dies off eventually.

{\it Organization:} Section \ref{sec:model} describes the model and Section \ref{sec:All-theorems} states the main results of the paper. It also derives the fluid limit model for the infection dynamics. Some of the proofs in the derivation of the fluid limit model are relegated to Appendix \ref{app:appendix-general-convergence}. Section \ref{sec:cost-equivalence-epsilon-0} reduces the fluid limit to a uni-dimensional dynamical system. Section \ref{sec:infection-sustains} identifies a sufficient condition on meeting rates and curing rates when infection sustains due to the presence of interactions with the mobile community. Section \ref{sec:opt-control-reduced-system} then derives the optimal control for the reduced dynamical system. Appendix \ref{app:equilibrium-point-concavity-convexity} contains some results used in Sections \ref{sec:cost-equivalence-epsilon-0} and \ref{sec:opt-control-reduced-system}. We end the paper with some concluding remarks in Section \ref{sec:conclusion}.

\section{Model}\label{sec:model}

\begin{figure*}
  \centering
  \includegraphics[scale=.5]{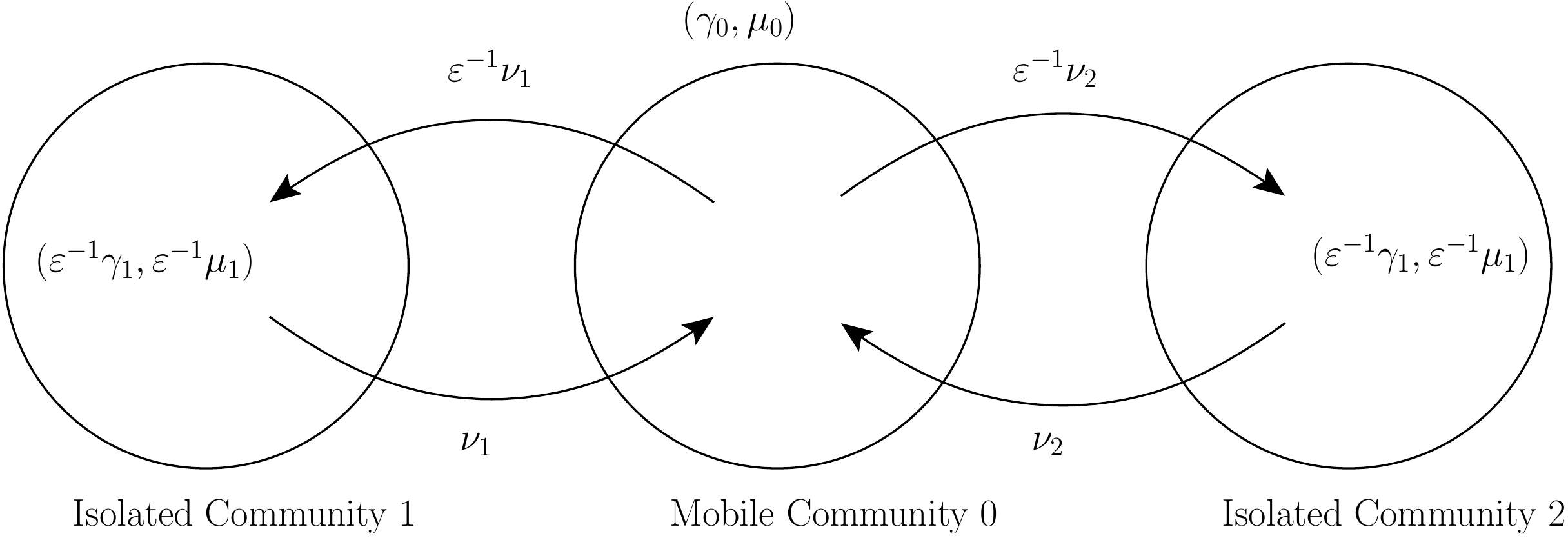}\\
  \caption{The interaction rates for a population model with two isolated communities and a mobile community.}
  \label{fig:interaction-rates}
\end{figure*}

Let us recall some notation from Section \ref{sec:intro}. We consider $K$ isolated communities and an influential mobile community. The members/nodes of each isolated community interact only among themselves and with members of the mobile community. There is no direct interaction with the members of another isolated community. We wish to study the spreading dynamics of an information/content among the members of the communities. A node that has a copy of the content is said to be infected and a node that does not have the copy of the content is said to be susceptible. When two nodes come in contact, one infected and the other susceptible, the former transmits a copy of the content to the latter, and the latter gets infected. An infected node may spontaneously get rid of the content via {\em curing}. The members of the isolated community interact within themselves at a \emph{faster} rate compared to the interactions within the members of the mobile community. This timescale separation leads to a nontrivial metastable phenomenon which we study in this paper.

We shall have two related scaling parameters $n$ and $\varepsilon$. The parameter $n$ will refer to population size (sum of populations of all the communities), which we shall drive to $\infty$ to arrive at a fluid limit. The parameter $\varepsilon$ will refer to timescale separation in the dynamics of interactions within the isolated communities and within the mobile community. The smaller the $\varepsilon$, the greater the separation. We shall drive $\varepsilon$ down to $0$. We shall also consider $\varepsilon$ as a function of $n$.

Let the $K$ isolated communities be indexed by $k=1,2,\ldots,K$. Refer to Figure \ref{fig:interaction-rates} that shows the case of two isolated communities and a mobile community. Community $k$ has size $M_k^n$. Any two members of community $k$ come in contact with each other at time instants that are points of an independent Poisson point process of rate $\varepsilon^{-1}\gamma_k^n$. An infected node in community $k$ spontaneously gets rid of the content at a rate $\varepsilon^{-1} \mu_k^n$. There is no direct contact between nodes of different communities for $k=1,2,\ldots,K$. However, there is \emph{indirect} inter-community contact via contact with the set of mobile nodes, which we shall call as community 0. This community is of size $M_0^n$. A mobile node has contacts with each node of community $k$ at instants that form an independent Poisson point process of rate $\varepsilon^{-1} \nu_k^n, 1 \leq k \leq K$. The spread of infection is {\em asymmetric}. Infection spread from mobile community $0$ to community $k$ happens at the indicated rate, but infection spread from community $k$ to mobile community $0$ however happens at a thinned rate of $\nu_k^n$ ($= \varepsilon \cdot (\varepsilon^{-1}\nu_k^n)$, hence ``thinned''). Any two mobile nodes come in contact with each other at instants that form a Poisson process of rate $\gamma^n_{0}$. An infected mobile node gets spontaneously cured at rate $\mu_0^n$.

Since $\varepsilon$ is going to be small, each isolated community with index $k=1,2,\ldots,K$ interacts within itself at a much faster timescale, and gets cured at a similarly faster timescale. Each of the $K$ communities also has a greater propensity to get infected via interactions with the mobile community. This fast timescale of interaction is visible from the dependence on $\varepsilon$ in the meeting rates $\varepsilon^{-1} \gamma_k^n$, $\varepsilon^{-1} \nu_k^n$, and in the curing rate $\varepsilon^{-1} \mu_k^n$, when $k = 1, \ldots, K$. The mobile community, however, is slow to be infected (within community contact rate is $\gamma_0^n$ and infection rate due to contact with community $k$ is $\nu_k^n$) and is slow to be cured (rate $\mu_0^n$). We will soon make some remarks on the asymmetry of this interaction.

We now allow the possibility of controlling the interactions of the mobile nodes by a campaigner. We do this by modulating the interaction rates $\nu_k^n,~k = 1,2,\ldots, K$ and $\gamma^n_0$ by a common \emph{control variable} $u(\cdot)$. The resulting interaction rates are $\varepsilon^{-1} \nu_k^n u(t)$ for $k = 1, \ldots, K$, and $\gamma_0^n  u(t)$. Let

\[\mathcal{U} = \{u:[0,T] \rightarrow [0,1]:u(\cdot) \mbox{ is measurable}\}.\]
An \emph{admissible} control policy must result in a control process with the following properties.
\begin{itemize}
 \item The control process sample path $u(\cdot) \in \mathcal{U}$ with probability $1$.
 \item The control process is non-anticipative (technically, $u(\cdot)$ is measurable with respect to the filtration ($\mathcal{F}_t,t\geq 0$) generated by the state and control variables up to the time that goes to define the filtration).
\end{itemize}

\begin{figure}
  \centering
  \includegraphics[scale=0.5]{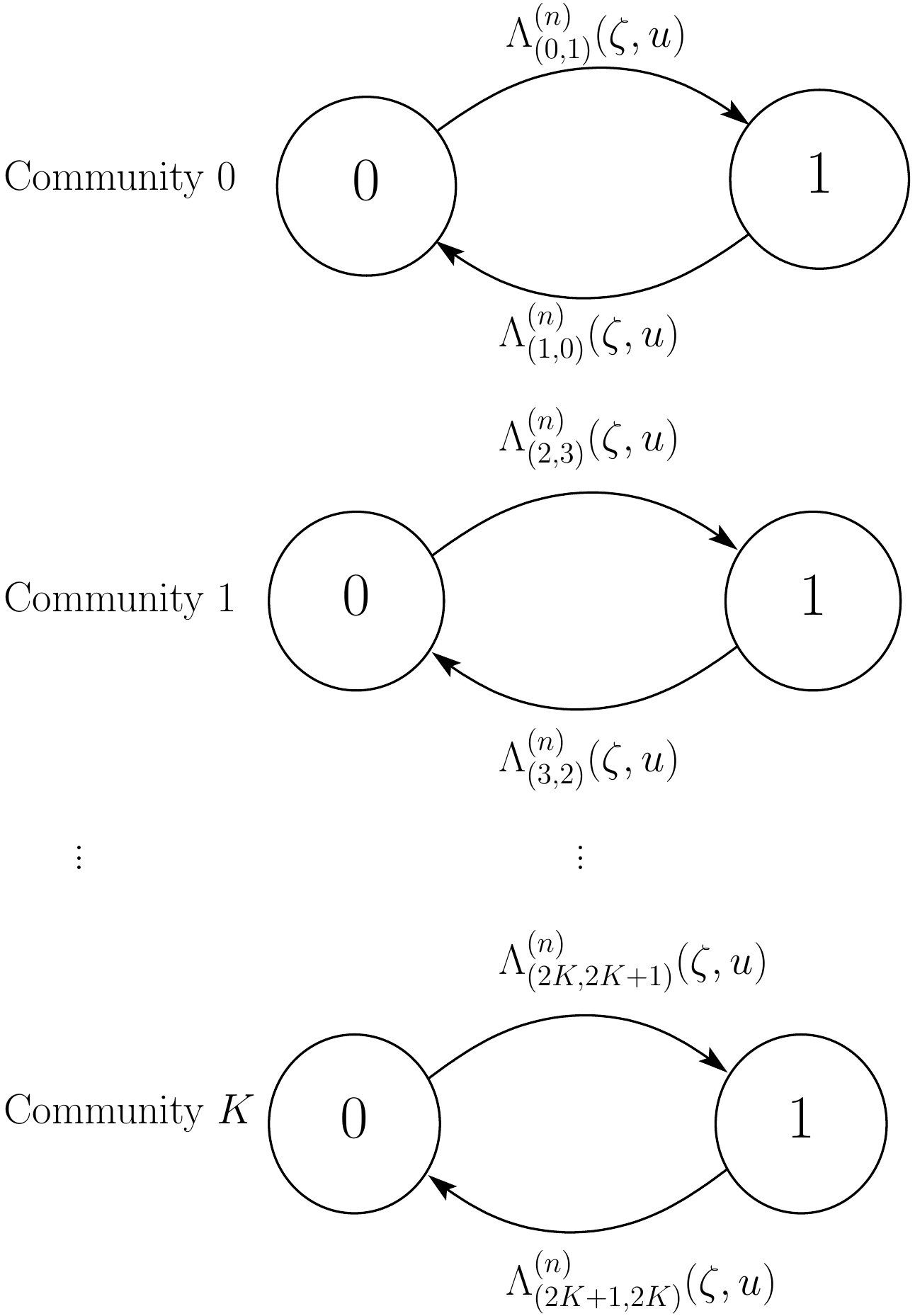}\\
  \caption{Allowed transition of a particle in the $(\varepsilon,n)$-system, along with the notation for transition rates.}\label{fig:transition-rates}
\end{figure}

We model the system of interacting communities as an $n$-particle closed system. Let the pair $(k,j)$ denote the state of a particle where $k$ indicates the community index of the particle and $j \in \{0,1\}$ indicates whether the particle is susceptible $(j=0)$ or infected $(j=1)$. See Figure \ref{fig:transition-rates}. Define $\mathcal{S} =\{(k,j):0\leq k \leq K,j \in \{0,1\}\}$. Let $Z^{n}_{p}(t)$ denote the state of the particle $p$ at time $t$; $\{Z^{n}_{p}(t),p=1,2,\ldots,n\}$ is a Markov decision process taking values in $\mathcal{S}$.

Let $\zeta^{n}(t)=(\zeta^n_{2k},\zeta^n_{2k+1}, 0 \leq k \leq K)$ denote the empirical distribution of the particles across states at time $t$, where $\zeta^{n}_{2k}(t)$ and $\zeta^{n}_{2k+1}(t)$ denote the fraction of susceptible nodes and infected nodes of the $k$th community, respectively. At time $t$, let $X_k^n(t)$ denote the number of nodes of the $k$th community that are infected. We then have
\begin{align}
\zeta^{n}_{2k}(t)&=\frac{M^{n}_{k}-X^{n}_{k}(t)}{n},0 \leq k \leq K,\nonumber\\
\zeta^{n}_{2k+1}(t)&=\frac{X^{n}_{k}(t)}{n}, 0 \leq k \leq K.
\end{align}
Let $\Lambda^{(n)}(\zeta,u) \in \mathbb{R}^{(2K+2) \times (2K+2)}$ denote the transition rate matrix of the $p$th particle's transitions, given the empirical distribution $\zeta^n(t)=\zeta$, defined as follows. Let $\Lambda^{(n)}_{(a,b)}(\zeta, u)$ denote the element in the $a$th row and $b$th column of the matrix $\Lambda^{(n)}(\zeta,u)$. Then $\Lambda^{(n)}_{(2k,2k+1)}(\zeta,u)$ is the rate of the transition from the state $(k,0)$ to $(k,1)$ and $\Lambda^{(n)}_{(2k+1,2k)}(\zeta,u)$ is the rate of the transition from the state $(k,1)$ to $(k,0)$. See Figure \ref{fig:transition-rates}. Our model for meeting rates and transitions is:
\begin{eqnarray}
\nonumber
\lefteqn{\Lambda^{(n)}_{(2k,2k+1)}(\zeta,u)}\\
\label{eqn:ratetranmatrix1}
& = & \left\{
\begin{array}{ll}
   n\zeta_{2k+1} \varepsilon^{-1}\gamma^{n}_{k}+n\zeta_{1}\varepsilon^{-1}\nu^{n}_{k}u ,& 1 \leq k \leq K,\\
       (n \zeta_{1}\gamma^{n}_{0}+\sum_{i=1}^{K}n\zeta_{2i+1}\nu^{n}_{i})u,& k=0,
\end{array}
\right.
\end{eqnarray}
and
\begin{equation}\label{eqn:ratetranmatrix2}
\Lambda^{(n)}_{(2k+1,2k)}(\zeta,u)=\left\{
\begin{array}{ll}
      \varepsilon^{-1}\mu^{n}_{k}, & 1 \leq k \leq K,\\
      \mu^{n}_{0} , & k=0.
\end{array}
\right.
\end{equation}
When $k=0$, we interpret (\ref{eqn:ratetranmatrix1}) as the meeting rate of an uninfected node in community $0$ with an infected node either in community $0$ or with one of the isolated communities. When $k \geq 1$, we interpret (\ref{eqn:ratetranmatrix1}) similarly. The parameter $u \in [0,1]$ is the campaigner's control. We interpret (\ref{eqn:ratetranmatrix2}) as spontaneous curing. No other transitions are allowed. This along with (\ref{eqn:ratetranmatrix1}) and (\ref{eqn:ratetranmatrix2}) completely define the transition rate matrix $\Lambda^{(n)}(\zeta,u)$.

We refer to the $n$-particle closed system with the timescale separation parameter $\varepsilon$ as the $(\varepsilon,n)$-system. Let $\mathcal{E}^{(n)}$ denote the space of empirical distributions of the $(\varepsilon,n)$-system:
\begin{eqnarray*}
\mathcal{E}^{(n)}\hspace*{-.1cm}=\{\zeta: \zeta_{2k},\zeta_{2k+1} \in \{0,1/n,\ldots,M^{n}_{k}/n\}, k=0,\ldots,K \\
 \text{ and } \zeta_{2k}+\zeta_{2k+1}=M^{n}_{k}/n, k=0,\ldots,K\}.
\end{eqnarray*}
The empirical distribution $\zeta^{n}(t)$ of the $(\varepsilon,n)$-system is also a continuous-time Markov decision process with state space $\mathcal{E}^{(n)}$ and control or action $u(t)$. Let $\zeta = (\zeta_0, \zeta_1, \ldots,\zeta_{2K},\zeta_{2K+1}) \in \mathcal{E}^{(n)}$ and let $\e^{n}_i$ be the vector with $1/n$ in the $i$th component and zeros elsewhere. The rates of transitions from $\zeta \mbox{ to } \zeta - \e^n_{2k}+\e^{n}_{2k+1}$ are easily seen to be
\begin{equation}
\begin{array}{ll}
  (M_k^n - n\zeta_{2k+1})(\varepsilon^{-1} \gamma_k^n n\zeta_{2k+1} + \varepsilon^{-1} \nu_k^n  n\zeta_1 u ), & k\not=0, \\
  (M_0^n - n\zeta_1) \left(\gamma_0^n n\zeta_1+\sum_{i=1}^K \nu_i^n n \zeta_{2i+1}\right)u, & k = 0.
\end{array}
\end{equation}
These denote increase in the number of infected nodes. The rates of transitions from state $\zeta$ to $\zeta + \e^{n}_{2k}- \e^{n}_{2k+1}$, which correspond to a decrease in the number of infected nodes, are given by
\begin{equation}
\begin{array}{ll}
  n\zeta_{2k+1} \varepsilon^{-1} \mu_k^n , & k\not=0, \\
  n\zeta_1 \mu_0^n , & k = 0.
  \end{array}
\end{equation}

For each $i \in \mS, j \in \mS$, the (right-continuous with left-limit) process $D_{i,j}^n(t)$ counts the number of $i$ to $j$ transitions in the time interval $[0,t]$. The evolution of the empirical distribution $\zeta^{n}(t)$ in the $(\varepsilon, n)$-system can be written as
\begin{equation}\label{eqn:empdisofnthsubsystem}
  \zeta^{n}(t) = \zeta^{n}(0) + \int_{[0,T]} \hspace*{-.1cm}\left[\Lambda^{(n)} \left( \zeta^{n}(s), u(s) \right)\right]^* \zeta^{n}(s) ~ds + \mathcal{M}^{n}(t),
\end{equation}
where $\mathcal{M}^{n}$ is a vector-valued square-integrable $(\F_t)$-measurable martingale whose $i$th component is
\begin{eqnarray*}
 \mathcal{M}^n_i(t) & = & \frac{1}{n} \sum_{j : j \neq i} D_{j,i}^n(t) - \frac{1}{n} \sum_{j': j' \neq i} D_{i,j'}^n(t) \\
 & & - \int_0^t \left[ \left[\Lambda^{(n)}(\zeta^{n}(s), u(s))\right]^* \zeta^{n}(s) \right]_i ~ds,
\end{eqnarray*}
and $\left[\Lambda^{(n)}(\cdot,\cdot)\right]^*$ denotes the adjoint of $\Lambda^{(n)}(\cdot,\cdot)$.

For an admissible control policy $\pi$, let $P^{\pi}$ be the induced probability measure associated with process of empirical measure, $\zeta^{n}(\cdot)$. The associated expectations are denoted $\mathbb{E}^{\pi}$. For a finite time horizon $T$, a controller employing policy $\pi$ incurs the following cost:
\begin{equation}
  \label{eqn:value-function}
  V^{(\varepsilon,n)}(\zeta^{n}(0), \pi) =  \mathbb{E}^{\pi} \Big[ \int_{[0,T]} \hspace*{-.15cm}r_1(\zeta^{n}(s), u(s)) ~ds + r_2(\zeta^{n}(T)) \Big],
\end{equation}
where $u(\cdot)$ is the control process resulting from the application of control policy $\pi$. The integrand $r_1(\zeta^{n}(s),u(s))$ constitutes the instantaneous cost and $r_2(\zeta^{n}(T))$ denotes the terminal cost for putting the empirical measure in state $\zeta^{n}(T)$ at terminal time $T$.
\section{Statements of the Main Results}\label{sec:All-theorems}
We are interested in, for well-separated timescales of small but nonzero $\varepsilon$, the evolution of the limiting system as the number of nodes increases. In the context of our model, this is a mean-field fluid limit. We shall assume the following on the population sizes:
\[
  n^{-1} M_k^n \rightarrow m_k.
\]
For convenience, we shall assume the following interaction rates:
\begin{align}
\begin{split}\label{eqn:asymptotic-rates}
  & n\gamma_k^n = \gamma_k,  k=0,1,\ldots,K \\
  & n \nu_k^n =\nu_k,  k= 1,2,\ldots,K. \\
  & \mu_k^n = \mu_k ,  k = 0,1,\ldots,K.
  \end{split}
\end{align}
See Appendix \ref{app:appendix-general-convergence} for a generalization where we only assume $n\gamma_k^n \rightarrow \gamma_k$, etc.
The appearance of $n$ in the {\em interaction} rates reflects the assumption that, given the $O(n)$ number of nodes in a community or a cross-community for interaction, the rate of increase of the fraction of susceptible nodes is $O(1)$. This is the regime where interesting interactions are visible in $O(1)$ time, if $\varepsilon$ is held constant while $n \rightarrow \infty$.

A natural fluid limit of the system in (\ref{eqn:empdisofnthsubsystem}) would be
\begin{equation}
  \label{eqn:fluidlimit}
  \zeta(t) = \zeta(0) + \int_{[0,t]}\left[ \Lambda \left( \zeta(s), u(s) \right)\right]^* \zeta(s) ~ds,
\end{equation}
where the matrix $\Lambda$ is obtained by taking limits in (\ref{eqn:ratetranmatrix1}) and (\ref{eqn:ratetranmatrix2}) and by using (\ref{eqn:asymptotic-rates}):
\begin{eqnarray}
\nonumber
\lefteqn{\Lambda_{(2k,2k+1)}(\zeta,u)}\\
\label{eqn:limitratetranmatrix1}
& = & \left\{
\begin{array}{ll}
    \left(\gamma_{0}\zeta_{1}+\sum_{i=1}^{K}\nu_{i}\zeta_{2i+1}\right)u(t),& k=0, \\
   \varepsilon^{-1}\gamma_{k}\zeta_{2k+1}+\varepsilon^{-1}\nu_{k}\zeta_{1}u(t) ,& 1 \leq k \leq K,
\end{array}
\right.
\end{eqnarray} and
\begin{equation}\label{eqn:limitratetranmatrix2}
\Lambda_{(2k+1,2k)}(\zeta,u)=\left\{
\begin{array}{ll}
  \mu_{0} , & k=0, \\
    \varepsilon^{-1}\mu_{k}, & 1 \leq k \leq K.
\end{array}
\right.
\end{equation}
We refer to the fluid limit as $(\varepsilon, \infty)$-system. Let $\mathcal{E}$ denote the space of empirical distributions of $(\varepsilon, \infty)$-system. The cost of employing a control policy $\pi$ for the $(\varepsilon, \infty)$-system (with a possibly random control process sample path) is
\[
  V^{(\varepsilon,\infty)}(\zeta(0), \pi) =  \mathbb{E}^{\pi} \Big[ \int_{[0,T]} \hspace*{-.15cm}r_1(\zeta(s), u(s)) ~ds + r_2(\zeta(T)) \Big],
\]
where the evolution of $\zeta$ is according to (\ref{eqn:fluidlimit}) for an initial condition $\zeta(0)$ and $u(\cdot)$ is the control sample path under the policy $\pi$. We assume that $r_1(\zeta,u)$ is decreasing in $\zeta$ (component-wise) and increasing in $u$ and that $r_2(\zeta)$ is decreasing in $\zeta$ (component-wise). Thus minimizing the cost function will lead to a combination of maximizing information diffusion at a fixed terminal time and minimizing control costs. The quantities $r_1$ and $r_2$ are chosen by the campaigner. They determine the tradeoff between infection spread and control costs. We make the following assumptions on the cost functions $r_1$ and $r_2$.

\vspace*{.2cm}
\begin{itemize}
 \item[({\bf A1})] The mapping $\zeta \mapsto r_1(\zeta, u)$ is Lipschitz continuous in $\zeta$, uniformly in $u \in [0,1]$. Similarly, the terminal reward mapping $\zeta \rightarrow r_2(\zeta)$ is also Lipschitz continuous. More precisely, there is an $L_1$ such that for every $\zeta, \zeta'\in \mathcal{E}$ and $u \in [0,1]$, we have
\begin{eqnarray*}
    |r_1(\zeta', u) - r_1(\zeta, u)| & \leq & L_1 ||\zeta - \zeta'||, \\
    |r_2(\zeta') - r_2(\zeta)| & \leq & L_1 || \zeta - \zeta'||.
  \end{eqnarray*}
\end{itemize}
We make the following assumptions on the parameters of the dynamics.

\vspace*{.2cm}
\begin{itemize}
\item[({\bf A2})]\label{ass:assumption-on-rates} $0 < m_k \gamma_k < \mu_k,~k=0,1,\ldots,K$.
\end{itemize}

\vspace*{.2cm}
We now relate the evolution of empirical distribution in the $(\varepsilon,n)$-system (\ref{eqn:empdisofnthsubsystem}) to that of the $(\varepsilon,\infty)$-system (\ref{eqn:fluidlimit}) when the same control sample path is applied to both the systems under the same initial condition $\zeta^{(\varepsilon,n)}(0) =\zeta^{(\varepsilon,\infty)}(0)=\zeta(0)$. In particular, we compare $\zeta^{(\varepsilon,n)}$ with $\zeta^{(\varepsilon,\infty)}$ and the corresponding costs.

\begin{theorem}\label{thm:cost-equality-n-infty}
Let $\zeta^{(\varepsilon,n)}(\cdot)$ denote the empirical distribution flow of the $(\varepsilon,n)$-system in (\ref{eqn:empdisofnthsubsystem}) when an admissible control policy $\pi_n$ is applied. Let $u(\cdot)$ be the corresponding control process sample path. Let  $\zeta^{(\varepsilon,\infty)}(\cdot)$ be the empirical  distribution flow of the  $(\varepsilon,\infty)$-system (\ref{eqn:fluidlimit}) resulting from the application of this (possibly random) control process sample path $u(\cdot)$ on the $(\varepsilon,\infty)$-system (\ref{eqn:fluidlimit}). Call the resulting control policy as $\pi_n$ again but on the $(\varepsilon,\infty)$-system. Assume $\zeta^{(\varepsilon,n)}(0)=\zeta^{(\varepsilon,\infty)}(0)=\zeta(0)$. Let $c > 0$. There are finite constants $0 < C, \bar{C} < \infty$ such that if $C/\log n \leq \varepsilon \rightarrow 0$ then
\begin{align}\label{eqn:theorem-p-bound}
 & P^{\pi_n} \left\{ \left\| \zeta^{(\varepsilon,n)}-\zeta^{(\varepsilon,\infty)} \right\|_{T} > \frac{c}{\log n} \right\} \leq \frac{\bar{C}}{\log n},  \\
 \label{eqn:theorem-v-bound}
 &\left \vert  V^{(\varepsilon,n)}(\zeta(0), \pi_n) - V^{(\varepsilon,\infty)}(\zeta(0), \pi_n) \right  \vert \leq \frac{\bar{C}}{\log n}.
\end{align}
\end{theorem}
\begin{IEEEproof}
See Appendix \ref{app:appendix-general-convergence}.
\end{IEEEproof}

\vspace*{.2cm}
The differential form of the $(\varepsilon,\infty)$-system in (\ref{eqn:fluidlimit}) is
\begin{align}\label{eqn:diff-form}
&\dot{\zeta}(t)=\left[\Lambda(\zeta(t),u(t))\right]^*\zeta(t), t \geq 0.
\end{align}
where $\dot{\zeta}(t)$ is the time derivative of $\zeta$. Since
\begin{equation}\label{eqn:commun-pop-constant}
\zeta_{2k}+\zeta_{2k+1}(t)= m_k, ~k=0,1,2,\ldots,K,
\end{equation} we have
\[\dot{\zeta}_{2k+1}(t)=-\dot{\zeta}_{2k}(t).\]
Hence it suffices to consider the evolution of components $\zeta_{2k+1}(t),k=0,1,2,\ldots,K$. Expanding (\ref{eqn:diff-form}) and using (\ref{eqn:commun-pop-constant}), the evolution of the empirical distribution of $(\varepsilon,\infty)$-system is the following:
\begin{align}
  \nonumber
  \varepsilon \dot{\zeta}_{2k+1}(t) = & - \mu_k \zeta_{2k+1}(t) \\
  \nonumber
  &  + (m_k - \zeta_{2k+1}(t)) \left( \gamma_k \zeta_{2k+1}(t) + \nu_k \zeta_{1}(t) u(t) \right), \\
  \label{eqn:drivingfunction-fast}
  &  \hspace*{3.8cm} k = 1, \ldots, K, \\
  \nonumber
  \dot{\zeta}_1(t)  = & - \mu_0 \zeta_1(t) \\
  \label{eqn:drivingfunction-slow}
  & + (m_0 - \zeta_1(t)) \Big( \gamma_0 \zeta_1(t) + \sum_{k=1}^K \nu_k \zeta_{2k+1}(t) \Big) u(t),\\
  \zeta_{2k}(t) = & ~ m_k -\zeta_{2k+1}(t),\quad k = 0,1, \ldots, K.\label{eqn:commun-pop-constant-2}
\end{align}
As $\varepsilon \rightarrow 0$, the fast timescale variables $\zeta_{2k+1}, 1 \leq k \leq K$ see the slow timescale variable $\zeta_{1}$ as a constant and rapidly converge to the equilibrium point of the dynamics in (\ref{eqn:drivingfunction-fast}) for a fixed $\zeta_1$. The limiting system as $\varepsilon \rightarrow 0$, denoted as $(0,\infty)$-system, is given by the following set of equations:
 \begin{align}
 \nonumber
 0 = & - \mu_k \zeta_{2k+1}(t) \\
 \nonumber
 & + (m_k - \zeta_{2k+1}(t)) \left( \gamma_k \zeta_{2k+1}(t) + \nu_k \zeta_{1}(t) u(t) \right),\\
 \label{eqn:limiting-epsilon-fast}
 & \hspace*{3.8cm}  k = 1, \ldots, K, \\
 \nonumber
  \dot{\zeta}_1(t) = & - \mu_0 \zeta_1(t) \\
  \label{eqn:limiting-epsilon-slow}
  & + (m_0 - \zeta_1(t)) \Big( \gamma_0 \zeta_1(t) + \sum_{k=1}^K \nu_k \zeta_{2k+1}(t) \Big) u(t),\\
  \zeta_{2k}(t) = & ~ m_k-\zeta_{2k+1}(t),\quad k = 0,1, \ldots, K.\label{eqn:commun-pop-constant-3}
\end{align}
Consider a specific case of instantaneous and terminal cost functions,
\begin{equation}
\label{eqn:cost-function}
\left.
\begin{array}{l}
r_1(\zeta(s),u(s))=u(s),  \\
r_2(\zeta(T))=-\zeta_{1}(T).
\end{array}
\right\}
\end{equation}
The remaining statements are specific to these cost functions\footnote{Other cost functions may also be of interest, for example, $-\int_0^T \zeta_1(s)~ds$ instead of $-\zeta_1(T)$. We restrict attention to the specified cost functions mainly for illustration of the key ideas.}. The quantity $r_1(\cdot,u)=u$ is the cost of a campaign of intensity $u$.
The cost of applying a control policy $\pi$ on the $(0,\infty)$-system (with a possibly random control process sample path) is
\begin{align}
V^{(0,\infty)}(\zeta(0), \pi)= \mathbb{E}^{\pi} \left[ \int_{[0,T]}  u(s) ~ds  - \zeta_{1}(T) \right]
\end{align}
where the evolution of $\zeta$ is according to (\ref{eqn:limiting-epsilon-fast})-(\ref{eqn:commun-pop-constant-3}) for an initial condition $\zeta(0)$ and $u(\cdot)$ is the control sample path under the policy $\pi$.

Let
\[\mathcal{U}_B=\{u(\cdot)\in \mathcal{U}: u(\cdot)\mbox{ has at most }B \mbox{ discontinuities}\}.\]
The remaining statements in this section consider the admissible control policies that result in a control process sample path $u(\cdot) \in \mathcal{U}_B$ with probability $1$.

We now relate the cost functions $V^{(\varepsilon,\infty)}$ and $V^{(0,\infty)}$ when identical control sample path is applied to the $(\varepsilon,\infty)$-system and the $(0,\infty)$-system.

\vspace*{.2cm}

\begin{theorem}\label{thm:cost-equality-epsilon-0}
Assume cost functions as in (\ref{eqn:cost-function}). Let $\zeta^{(\varepsilon,\infty)}(\cdot)$ denote the empirical distribution flow of the $(\varepsilon,\infty)$-system in (\ref{eqn:drivingfunction-fast})-(\ref{eqn:commun-pop-constant-2}) when an admissible control policy $\pi$ is applied. Let $u(\cdot)$ be the corresponding (possibly random) control process sample path. Let $u(\cdot)$ have at most $B$ discontinuities, i.e., $u(\cdot) \in \mathcal{U}_B$, with probability $1$. Let $\zeta^{(0,\infty)}(\cdot)$ be the empirical distribution flow of the  $(0,\infty)$-system (\ref{eqn:limiting-epsilon-fast})-(\ref{eqn:commun-pop-constant-3}) resulting from the application of the same control process sample path on the $(0,\infty)$-system. Call this policy as $\pi$ again but on the $(0,\infty)$-system. Assume $\zeta^{(\varepsilon,\infty)}(0)=\zeta^{(0,\infty)}(0)=\zeta(0)$. Then
\begin{align}
&\sup_{t \in [0,T]}\left\vert  \zeta_1^{(\varepsilon,\infty)}(t)-\zeta_1^{(0,\infty)}(t)\right\vert \leq C_1(B) \varepsilon \label{eqn:cost-equality-epsilon-0-2} \mbox{ with probability } 1,
\end{align}
where $C_1(B)$ is a system dependent constant that depends on the maximum number of discontinuities of $u(\cdot)$, but does not depend on $\varepsilon$. We also have
\begin{align}
&\vert  V^{(\varepsilon,\infty)}(\zeta(0), \pi)-V^{(0,\infty)}(\zeta(0), \pi) \vert \leq C_1(B) \varepsilon.\label{eqn:cost-equality-epsilon-0-1}
\end{align}
\end{theorem}
\begin{IEEEproof}For a proof, see Section \ref{sec:proof-via-kokotovic}.
\end{IEEEproof}
\vspace*{.2cm}
%

As a consequence of Theorems \ref{thm:cost-equality-n-infty} and \ref{thm:cost-equality-epsilon-0}, we have the following result.

\begin{theorem}\label{thm:cost-equality-epsilon-n-0-infty}
Let $c > 0$. Let $C$ and $\bar{C}$ be as in Theorem \ref{thm:cost-equality-n-infty}, let $C_1(B)$ be as in Theorem \ref{thm:cost-equality-epsilon-0}, and set $\varepsilon = C/(\log n)$. Assume cost functions as in (\ref{eqn:cost-function}). Let $\zeta^{(\varepsilon,n)}(\cdot)$ denote the empirical distribution flow of the $(\varepsilon,n)$-system (\ref{eqn:empdisofnthsubsystem}) when an admissible control policy $\pi$ is applied. Let $u(\cdot)$ be the corresponding control process sample path. Let $u(\cdot)$ have at most $B$ discontinuities, i.e., $u(\cdot) \in \mathcal{U}_B$, with probability $1$. Let  $\zeta^{(0,\infty)}(\cdot)$ be the empirical  distribution flow of the  $(0,\infty)$-system (\ref{eqn:limiting-epsilon-fast})-(\ref{eqn:commun-pop-constant-3}) resulting from the application of the same control process sample path $u(\cdot)$ on the $(0,\infty)$-system. Call the resulting policy as $\pi$ again but on the $(0,\infty)$-system. Assume $\zeta^{(\varepsilon,n)}(0)=\zeta^{(0,\infty)}(0)=\zeta(0)$. Then
\begin{align}\label{eqn:infection-equality-epsilon-n-0-infty}
&P^{\pi} \left\{ \left\| \zeta_1^{(\varepsilon,n)}-\zeta_1^{(0,\infty)} \right\|_{T} > \frac{c + C_1(B)C}{\log n} \right\} \leq \frac{\bar{C}}{\log n},  \\
\label{eqn:cost-equality-epsilon-n-0-infty}
&\left \vert  V^{(\varepsilon,n)}(\zeta(0), \pi) - V^{(0,\infty)}(\zeta(0), \pi) \right \vert \leq \frac{\bar{C} + C_1(B)C}{\log n}.
\end{align}
\end{theorem}
\begin{IEEEproof}
The proof is straightforward. Consider a policy $\pi$, apply Theorem \ref{thm:cost-equality-epsilon-0}, and apply Theorem \ref{thm:cost-equality-n-infty} with $\pi_n = \pi$. We then have
\begin{eqnarray*}
\lefteqn{ P^{\pi} \left\{ \left\| \zeta_1^{(\varepsilon,n)}-\zeta_1^{(0,\infty)} \right\|_{T} > \frac{c + C_1(B)C}{\log n} \right\} } \\
& \leq & P^{\pi} \left\{ \left\| \zeta_1^{(\varepsilon,n)}-\zeta_1^{(\varepsilon,\infty)} \right\|_{T} +\left\| \zeta_1^{(\varepsilon,\infty)}-\zeta_1^{(0,\infty)} \right\|_{T} \right.\\
& & \left.  \hspace*{4.35cm} > \frac{c + C_1(B)C}{\log n} \right\} \\
& \leq & P^{\pi} \left\{ \left\| \zeta_1^{(\varepsilon,n)}-\zeta_1^{(\varepsilon,\infty)} \right\|_{T} >\frac{c}{\log n} \right.\\
& & \left. \hspace*{0.7cm} \mbox{ or }\left\| \zeta_1^{(\varepsilon,\infty)}-\zeta_1^{(0,\infty)} \right\|_{T} > \frac{C_1(B)C}{\log n} \right\} \\
& \leq & P^{\pi} \left\{ \left\| \zeta_1^{(\varepsilon,n)}-\zeta_1^{(\varepsilon,\infty)} \right\|_{T} >\frac{c}{\log n} \right\} \\
& & + P^{\pi} \left\{ \left\| \zeta_1^{(\varepsilon,\infty)}-\zeta_1^{(0,\infty)} \right\|_{T} > C_1(B) \varepsilon \right\} \\
& \leq & \frac{\bar{C}}{\log n} + 0\\
& = & \frac{\bar{C}}{\log n},
\end{eqnarray*}
where the penultimate inequality follows from (\ref{eqn:theorem-p-bound}) and (\ref{eqn:cost-equality-epsilon-0-2}). This proves (\ref{eqn:infection-equality-epsilon-n-0-infty}).

To get (\ref{eqn:cost-equality-epsilon-n-0-infty}), add and subtract the term $V^{(\varepsilon,\infty)}(\zeta(0), \pi)$ inside the modulus of the left-hand side  of (\ref{eqn:cost-equality-epsilon-n-0-infty}), apply the triangle inequality, and use equations (\ref{eqn:theorem-v-bound}) and (\ref{eqn:cost-equality-epsilon-0-1}) to bound the terms.
\end{IEEEproof}

\vspace*{.2cm}

The following theorem gives a characterization of the optimal control for the $(0,\infty)$-system.

\vspace*{.2cm}

\begin{theorem}\label{thm:opt-control}
Assume cost functions as in (\ref{eqn:cost-function}). On the $(0,\infty)$-system, the  open-loop deterministic policy
\begin{equation}\label{eqn:optimal-control}
u^{\star}(t)=\left\{
\begin{array}{ll}
  0, & 0 \leq t < \tau, \\
    1, & \tau \leq t \leq T,
\end{array}\right.
\end{equation}
denoted $\hat{\pi}_{(0,\infty)}$, is optimal for some $\tau \in [0,T]$.
\end{theorem}

\begin{IEEEproof} See Section \ref{sec:opt-control-reduced-system}.
\end{IEEEproof}

\vspace*{.2cm}

Finally, we show that the cost incurred in applying the control $u^{\star}(\cdot)$ on the $(\varepsilon,n)$-system is within $O(1/ \log n)$ of the optimal cost of $(\varepsilon,n)$-system. Let
\begin{align}\label{eqn:optimal-control-epsilon-n-B}
V_{\star,B}^{(\varepsilon,n)}=\inf_{\substack{\pi: \pi \text{ is admissible, }\\ u(\cdot) \in \mathcal{U}_B \text{ w.p.}1}}V^{(\varepsilon,n)}(\zeta(0), \pi )
\end{align}

\vspace*{.2cm}

\begin{theorem}
Let $B \geq 1$. Assume cost functions as in (\ref{eqn:cost-function}). Let $C$ and $\bar{C}$ be as in Theorem \ref{thm:cost-equality-n-infty}, let $C_1(B)$ be as in Theorem \ref{thm:cost-equality-epsilon-0}, and set $\varepsilon = C/(\log n)$. Let $\hat{C} = \bar{C} + C_1(B)C$. Let $\hat{\pi}_{(0,\infty)}$ be the optimal control for the $(0,\infty)$-system  with $u^{\star}(\cdot)$ of Theorem \ref{thm:opt-control} as given in (\ref{eqn:optimal-control}) with $\tau \in [0,T]$. Then
\begin{align}\label{eqn:mainresult}
\left|
V_{\star,B}^{(\varepsilon,n)} - V^{(\varepsilon,n)}(\zeta(0), \hat{\pi}_{(0,\infty)})
\right| \leq \frac{3\hat{C}}{\log n}.
\end{align}
\end{theorem}
\begin{IEEEproof}
This follows by applying Theorem \ref{thm:cost-equality-epsilon-n-0-infty} twice as follows. Let $\hat{\pi}_{(\varepsilon,n)}$ be a control for the $(\varepsilon,n)$-system resulting in a control process sample path $u(\cdot) \in \mathcal{U}_B$ with probability 1 and satisfying
\begin{align}\label{eqn:near-optimal-control-epsilon-n-B}
\left \vert V_{\star,B}^{(\varepsilon,n)}-V^{(\varepsilon,n)}(\zeta(0), \hat{\pi}_{(\varepsilon,n)})\right \vert \leq \frac{\hat{C}}{\log n}.
\end{align}
We then have
\begin{align*}
V^{(\varepsilon,n)}(\zeta(0), \hat{\pi}_{(0,\infty)}) &\stackrel{(a)}{\geq} V_{\star,B}^{(\varepsilon,n)}\\
& \stackrel{(b)}{\geq} V^{(\varepsilon,n)}(\zeta(0), \hat{\pi}_{(\varepsilon,n)}) - \frac{\hat{C}}{\log n}\\
 &\stackrel{(c)}{\geq} V^{(0,\infty)}(\zeta(0), \hat{\pi}_{(\varepsilon,n)}) - \frac{2\hat{C}}{\log n}\\
  &\stackrel{(d)}{\geq} V^{(0,\infty)}(\zeta(0), \hat{\pi}_{(0,\infty)}) - \frac{2\hat{C}}{\log n}\\
  &\stackrel{(e)}{\geq} V^{(\varepsilon,n)}(\zeta(0), \hat{\pi}_{(0,\infty)}) - \frac{3\hat{C}}{\log n}.
\end{align*}
The inequality (a) is due to (\ref{eqn:optimal-control-epsilon-n-B}) and the fact that $u^{\star}(\cdot) \in \mathcal{U}_B$. The inequality (b) is due to (\ref{eqn:near-optimal-control-epsilon-n-B}). The inequality (c) is due to Theorem \ref{thm:cost-equality-epsilon-n-0-infty} with $\pi=\hat{\pi}_{(\varepsilon,n)}$. The inequality (d) is due to the optimality of the control $\hat{\pi}_{(0,\infty)}$ on the $(0,\infty)$-system. The inequality (e) is due to Theorem \ref{thm:cost-equality-epsilon-n-0-infty} this time applied with $\pi=\hat{\pi}_{(0,\infty)}$. Summarizing, we have
\begin{align*}
V^{(\varepsilon,n)}(\zeta(0), \hat{\pi}_{(0,\infty)})
& \geq V_{\star,B}^{(\varepsilon,n)} \\
& \geq V^{(\varepsilon,n)}(\zeta(0), \hat{\pi}_{(0,\infty)}) - \frac{3\hat{C}}{\log n}.
\end{align*}
This completes the proof.
\end{IEEEproof}

\vspace*{.2cm}

The upshot is that an open loop threshold policy $u^{\star}(\cdot)$, which is optimal on the $(0,\infty)$-system, is asymptotically optimal, among policies that yield control process sample paths in $\mathcal{U}_B$ w.p. 1, on the $(\varepsilon=C/(\log n),n)$-system with an error of the order $O(1/\log n)$.

\section{Cost Equivalence of $(\varepsilon,\infty)$ and $(0,\infty)$ Systems}
\label{sec:cost-equivalence-epsilon-0}

\subsection{Preliminaries}
We first recall the $(\varepsilon,\infty)$-system of (\ref{eqn:drivingfunction-fast})-(\ref{eqn:commun-pop-constant-2}). As a result of (\ref{eqn:commun-pop-constant}), the empirical distribution $\zeta(t)$ of $(\varepsilon,\infty)$-system is completely defined by the components $\zeta_{2k+1}(t),k=0,1,\ldots,K$. Therefore, we consider the evolution of these components alone. Denote $x_k(t)=\zeta_{2k+1}(t),~ k=0,1,\ldots,K$, for convenience. The evolution of $x_k(t), k=0,1,\ldots,K$ according to (\ref{eqn:drivingfunction-fast}) and (\ref{eqn:drivingfunction-slow}) is
\begin{align}
  \nonumber
  \varepsilon \dot{x}_k(t)  = & - \mu_k x_k(t) \\
  \nonumber
  & + (m_k - x_k(t)) \cdot \left( \gamma_k x_k(t) + \nu_k x_0(t) u(t) \right), \\
  \label{eqn:drivingfunction-fast-1}
  & \hspace*{4.3cm} k = 1, \ldots, K,
\end{align}
\begin{align}
  \nonumber
  \dot{x}_0(t) = & - \mu_0 x_0(t) \\
  \label{eqn:drivingfunction-slow-1}
  & + (m_0 - x_0(t)) \Big( \gamma_0 x_0(t) + \sum_{k=1}^K \nu_k x_k(t) \Big) u(t).
\end{align}

Define $x(t) := ( x_k(t) )_{k = 1, \ldots, K}$, define
\begin{align}
  \nonumber
  g_k(x_0, x_k, u)  := & - \mu_k x_k + (m_k - x_k) \cdot \left( \gamma_k x_k + \nu_k x_0 u \right),\\
  \label{eqn:drivingfunction-fast-rhs}
  & \hspace*{3cm} k = 1, \ldots, K,
\end{align}
and further define
\begin{equation}
  \label{eqn:drivingfunction-fast-rhs-vector}
  g(x_0, x, u) := \left( g_k(x_0, x_k, u) \right)_{k = 1, \ldots, K}.
\end{equation}
Similarly, define
\[
  f(x_0, x, u) := - \mu_0 x_0 + (m_0 - x_0) \Big( \gamma_0 x_0 + \sum_{k=1}^K \nu_k x_k \Big) \cdot u.
\]
Then we can rewrite the $(\varepsilon,\infty)$-system in (\ref{eqn:drivingfunction-fast})-(\ref{eqn:drivingfunction-slow}), or equivalently in (\ref{eqn:drivingfunction-fast-1})-(\ref{eqn:drivingfunction-slow-1}), compactly as
\begin{align}
\label{eqn:epsilon-infty-compact-1}
  \varepsilon \dot{x}(t) & =  g(x_0(t), x(t), u(t)) \\
  \dot{x}_0(t) & =  f(x_0(t), x(t), u(t)). \label{eqn:epsilon-infty-compact-2}
\end{align}
Similarly, we can write the $(0,\infty)$-system in (\ref{eqn:limiting-epsilon-fast}) and (\ref{eqn:limiting-epsilon-slow}) as
\begin{align}
0 & =  g_k(x_0(t), x_k(t), u(t)),~k=1,\ldots,K, \label{eqn:equilibrium-point}\\
  \dot{x}_0(t) & =  f(x_0(t), x(t), u(t)).\nonumber
\end{align}
Recall the intuition that as $\varepsilon \rightarrow 0$, the fast timescale variables $x_k$ of the $(\varepsilon,\infty)$-system see the slow timescale variable $x_0$ as a constant and rapidly approach the equilibrium associated with this fixed $x_0$. This equilibrium, denoted as $x_k^\star$, is the $k^{th}$ component of the solution to (\ref{eqn:equilibrium-point}); it depends on $x_0$ only through $\xi_k = \nu_k x_0 u$, and so we write $x_k^\star(\xi_k)$. By factoring the quadratic form $g_k(x_0, x_k, u)$ in (\ref{eqn:drivingfunction-fast-rhs}), for $\gamma_k\neq 0$, we get
\begin{equation}
  \label{eqn:g_k-quadratic}
  g_k(x_0, x_k, u) = -\gamma_k (x_k - x_k^+(\xi_k)) (x_k + x_k^-(\xi_k))
\end{equation}
where, with $b_k := \mu_k - m_k \gamma_k$, we have
\begin{align}
  \label{eqn:xk*}
  x_k^+(\xi_k) = & ~\frac{-(b_k + \xi_k) + \sqrt{(b_k + \xi_k)^2 + 4 \gamma_k m_k  \xi_k}}{2 \gamma_k}, \\
  \label{eqn:xk*bar}
  x_k^-(\xi_k) = & ~\frac{(b_k + \xi_k) + \sqrt{(b_k + \xi_k)^2 + 4 \gamma_k m_k \xi_k }}{2 \gamma_k}.
\end{align}
By Assumption ({\bf A2}), $m_k \gamma_k < \mu_k$ and so $b_k > 0$. It is also clear from (\ref{eqn:xk*}) that $x_k^{+}(\xi_k)\geq 0$ and that equality holds if and only if $\xi_k =0$. Thus, $x_k^+(0) = 0$. Applying the inequality $\sqrt{1+y}\leq 1+y/2$, $y \geq 0$, to (\ref{eqn:xk*}), we have $x_k^+(\xi_k) \leq \frac{\xi_k}{b_k+\xi_k}m_k < m_k$. We have therefore verified that $0\leq x_k^+(\xi_k)< m_k.$

Since $x_k^-(\xi_k)>0$, there is exactly one solution to $g_k(x_0, \cdot, u) = 0$ in the interval $[0, m_k]$, and this is $x_k^+(\xi_k)$.

For each fixed $x_0$ and $u$, $g_k(x_0, x_k, u) >0$ for $0 \leq x_k < x_k^+(\xi_k)$ and $g_k(x_0, x_k, u) <0$ for $x_k^+(\xi_k) < x_k \leq m_k $. Hence, for each fixed $x_0$ and $u$, {\it the point $x_k^+(\xi_k)$ is the globally (in $[0,m_k]$) asymptotically stable equilibrium for the dynamics (\ref{eqn:epsilon-infty-compact-1})}.
Thus $x_k^{\star}(\xi_k)=x_k^{+}(\xi_k), ~k=1,2,\ldots,K$.

The following hold true for the equilibrium point $x_k^{\star}(\xi_k)$.
\begin{itemize}
 \item The mapping $\xi_k \mapsto x^\star_k(\xi_k)$ is a strictly increasing and strictly concave function. See Proposition \ref{prop:equilibrium-point-concavity} of Appendix \ref{app:equilibrium-point-concavity-convexity}.
  \item The mapping $\xi_k \mapsto \xi_k x^\star_k(\xi_k)$ is a strictly increasing, but now a strictly convex function. See Proposition \ref{prop:equilibrium-point-convexity} of Appendix \ref{app:equilibrium-point-concavity-convexity}.
\end{itemize}

Intuitively, then, as $\varepsilon \rightarrow 0$, the slow timescale variables see the fast variables as having equilibrated to $x_k^{\star}(\xi_k)=x_k^{\star}(\nu_k x_0 u)$. Thus, intuitively, the evolution of $(0,\infty)$-system is given by the reduced system:
\begin{align}\label{eqn:0-infty-equilibrium}
\dot{x}_0(t) & =  f(x_0(t), x^\star(\xi(t)), u(t)),
\end{align}
where
\begin{eqnarray*}
\xi(t)  & = & (\xi_k(t))_{k=1,2,\ldots,K}, \\
x^{\star}(\xi(t)) & = & (x_k^{\star}(\xi_k(t)))_{k=1,2,\ldots,K}, \\
\xi_k(t) & = & \nu_kx_0(t)u(t), k = 1, \ldots, K.
\end{eqnarray*}
With these we now rigorously prove Theorem \ref{thm:cost-equality-epsilon-0}.

\subsection{Proof of Theorem \ref{thm:cost-equality-epsilon-0}}\label{sec:proof-via-kokotovic}
We prove Theorem \ref{thm:cost-equality-epsilon-0} by appealing to Kokotovi\'{c}'s \cite[Thm. 2.1]{kokotovic1984applications}. The following hold true for the fast dynamics in (\ref{eqn:epsilon-infty-compact-1}).
\begin{itemize}
  \item For a fixed $x_0$ and $u$, the point $x_k^*(\xi_k)$ is a globally asymptotically stable equilibrium for the dynamics (\ref{eqn:epsilon-infty-compact-1}). Further, the asymptotic stability is uniform in $x_0$ and $u$. Global asymptotic stability follows from the remarks made earlier. That the asymptotic stability is uniform in $x_0$ and $u$ follows from the assumption ({\bf A2}) that $\gamma_k m_k < \mu_k$ for all $k$. 
  \item Let $x := (x_k)_{k=1,2,\ldots,K}$. The Jacobian matrix $\frac{\partial g}{\partial x}$ is diagonal, and its eigenvalues (which are indeed real) are all strictly on the left half plane. This follows from the fact that $g_k$ is quadratic and concave in $x_k$ for fixed $x_0,u$.
  \end{itemize}
These are the two assumptions needed to apply Kokotovi\'{c}'s \cite[Thm. 2.1]{kokotovic1984applications}. Let $x(\cdot):=(x_k(\cdot))_{k=0,1,\ldots,K}$ denote the solution to $(\varepsilon,\infty)$-system in (\ref{eqn:epsilon-infty-compact-1}) and (\ref{eqn:epsilon-infty-compact-2}) when an admissible control policy $\pi$ is applied, with initial condition $x_0(0)=\zeta_1(0),x_k(0)=\zeta_{2k+1}(0),~ k=1,2,\ldots,K$. Let $u(\cdot)$ be the resulting control sample path. Let $\overline{x}_0(\cdot)$ be the solution to $(0,\infty)$-system in (\ref{eqn:0-infty-equilibrium}) on applying $u(t)$.
Take $\overline{x}_0(0)=x_0(0)=\zeta_{1}(0)$ and $x_k(0) =\zeta_{2k+1}(0),~k=1,2,\ldots,K$. Then $x_0(t)=\zeta^{(\varepsilon,\infty)}_{1}(t)$ and $\overline{x}_0(t)=\zeta_{1}^{(0,\infty)}(t)$.
Kokotovi\'{c}'s \cite[Thm. 2.1]{kokotovic1984applications} then bounds the error between $x_0(\cdot)$ and $\overline{x}_0(\cdot)$ for each sample path as
\begin{align}
\label{eqn:kokotovic-bound}
\sup_{t \in [0,T]} \left \vert x_0(t)-\overline{x}_0(t)\right \vert  =O(\varepsilon).
\end{align}
We must now argue that the constant multiplier for $O(\varepsilon)$ in (\ref{eqn:kokotovic-bound}) is independent of the control process sample path. To do this, we exploit the finite number of discontinuities assumption and appeal to a result in \cite[Thm.~1.1]{levin1954singular} to get a  refinement of the above statement:
\begin{align}\label{eqn:error-bound-control-independent}
\sup_{t \in [0,T]} \left \vert x_0(t)-\overline{x}_0(t)\right \vert  \leq C_1(B) \varepsilon,
\end{align}
where $C_1(B)$ is independent of the control process sample path $u(\cdot)$ in $\mathcal{U}_B$.

To show (\ref{eqn:error-bound-control-independent}), we proceed as follows. Suppose that, in addition to the assumptions needed for Kokotovi\'{c}'s result \cite[Thm.~2.1]{kokotovic1984applications}, we also have that $f,g,\frac{\partial f}{\partial x_0},\frac{\partial f}{\partial x},\frac{\partial g}{\partial x_0},\frac{\partial g}{\partial x}$ are continuous in $(x_0,x,u,t)$. Then, from \cite[Eqn.~(3.24)]{levin1954singular} in the proof of \cite[Thm.~1.1]{levin1954singular} restated in our notation, we get
\begin{align}\label{eqn:error-bound-epsilon-0}
 \left \vert x_0(t)-\overline{x}_0(t)\right \vert \leq \left \vert x_0(t_0)-\overline{x}_0(t_0)\right \vert \exp\{C_2(t-t_0)\}+C_3 \varepsilon
\end{align}
for $t \in [t_0,t_0+l]$ for some $l > 0$. The constants $C_2,C_3 $ and $l$ are independent of $t_0$ and $\varepsilon$. We now apply (\ref{eqn:error-bound-epsilon-0}) repeatedly to each continuous segment of the control process sample path or a part thereof so that each segment is of length at most $l$. The boundedness of $u(\cdot)$ allows us to use the same $C_2$ and $C_3$ across all sample paths. The total number of line segments is
$$\sum_{i=1}^B \lceil l_i/l \rceil \leq \sum_{i=1}^B (l_i/l +1) = T/l + B,$$
where $l_i$ are the lengths of the $B$ continuous segments of the sample path $u(\cdot)$, some of which may be zero. Applying the bound (\ref{eqn:error-bound-epsilon-0}) to each of these segments and telescoping, we get
\begin{align}
\nonumber
\sup_{t \in [0,T]} \left \vert x_0(t)-\overline{x}_0(t) \right \vert
& \leq C_3 \varepsilon \cdot \frac{\exp\{(T/l+B)C_2 l\} - 1}{\exp\{C_2 l\}-1} \\
\label{eqn:error-bound-control-independent-1}
& = C_1(B) \varepsilon.
\end{align}
This establishes (\ref{eqn:error-bound-control-independent}).

As a consequence of (\ref{eqn:error-bound-control-independent}), for cost functions in (\ref{eqn:cost-function}), we have
\begin{eqnarray*}
\lefteqn{\left| V^{(\varepsilon,\infty)}(\zeta(0), \pi)-V^{(0,\infty)}(\zeta(0),\pi) \right| } \\
& = & \mathbb{E}^{\pi} \left[ \left \vert \zeta^{(\varepsilon,\infty)}_{1}(T)-\zeta^{(0,\infty)}_{1}(T)\right \vert \right] \\
& \leq & \mathbb{E}^{\pi} \left[ \sup_{t \in [0,T]} \| x_0(t)-\overline{x}_0(t) \| \right] \\
& \leq & C_1(B)\varepsilon.
\end{eqnarray*}
This completes the proof of Theorem \ref{thm:cost-equality-epsilon-0}.
$\hfill \IEEEQEDclosed$

\section{A sufficient condition for infection to sustain}
\label{sec:infection-sustains}
We now investigate the benefit of cross-community interactions in sustaining the infection. As a consequence of Theorem \ref{thm:cost-equality-epsilon-n-0-infty}, for models with sufficiently large population size $n$ and small timescale separation parameter $\varepsilon$, the infection level in the mobile community on the $(\varepsilon,n)$-system, $\zeta^{(\varepsilon,n)}_1$, and that on the $(0,\infty)$-system, $\zeta^{(0,\infty)}_1$, are close to each other with high probability. Hence the study of the evolution of $\zeta^{(0,\infty)}_1$ (or $x_0(t)$ in the simplified notation introduced in this section) will help us gain insight on the spread of infection in the original population model. The dynamics on the $(0,\infty)$-system, from (\ref{eqn:0-infty-equilibrium}), is
\begin{align}\label{eqn:0-infty-equilibrium-recall}
\dot{x}_0(t) & =  f(x_0(t), x^\star(\xi(t)), u(t)),
\end{align}
where
\begin{align}
  \nonumber
  f( &x_0, x^\star(\xi), u) \\
  \label{eqn:actual-dynamics}
  & = - \mu_0 x_0 + (m_0 - x_0) \Big( \gamma_0 x_0 + \sum_{k=1}^K \nu_k x^{\star}_k(\xi_k) \Big) \cdot u, \\
  \xi_k & = ~\nu_k x_0 u, \quad k = 1, \ldots, K.
\end{align}
Fix $u=1$. We consider the evolution of $x_0(t)$ in the following cases.
\begin{itemize}
\item[(a)] Cross-community interactions are absent, i.e., $\nu_k =0, k=1,\ldots,K$. Alternatively, $u(t) \equiv 0$. Observe that $\xi_k(t)=\nu_k x_0(t) u(t) \equiv 0$ and hence $x_k^{\star}(\xi_k(t)) \equiv 0, ~k=1,\ldots,K$. So $x^{\star}(\xi(t))={\bf 0}$, the all-zero vector. Then
\begin{align}
f(x_0, {\bf 0}, {\bf 0}) &= - \mu_0 x_0 + (m_0 - x_0) \gamma_0 x_0 \label{eqn:zero-cross-community-upperbound}\\
                        &\stackrel{(a)}{\leq} - \mu_0 x_0 + m_0 \gamma_0 x_0 \nonumber\\
                       & \stackrel{(b)}{\leq} 0,\nonumber%
\end{align}
where (b) is due to Assumption ({\bf A2}). Further, the inequalities (a) and (b) are satisfied with equality if and only if $x_0 = 0$. Therefore, $f(x_0, {\bf 0}, {\bf 0}) = 0$ for $x_0=0$ and $f(x_0, {\bf 0}, {\bf 0}) < 0$ for $x_0 > 0$. Thus $x_0 = 0$ is a stable equilibrium point, and the infection dies out for any initial condition.

 \item[(b)] Cross-community interaction rates $\nu_k >0,  k=1,2,\ldots, K$. Take $u(\cdot) \equiv 1$, also written as $u(\cdot) = {\bf 1}$. The equilibrium point $x_0$ satisfies
 \begin{align*}
  f(x_0, x^\star(\xi), {\bf 1}) = 0, \quad \xi = (\nu_k x_0)_{k=1,\ldots, K}.
 \end{align*}
Figure \ref{fig:network-effect-helps} plots $f(x_0, x^\star(\xi), {\bf 1})$ for the case of two isolated communities and a mobile community. The fraction of total population in each community is $m_0=0.4, m_1=m_2=0.3$. The within-community interaction rate parameters are $\gamma_0=\gamma_1=\gamma_2=1$. The curing rate parameters are $\mu_0=\mu_1=\mu_2=2$. The solid line shows the case when cross-community interactions are present with rate parameters $\nu_1=\nu_2=8$.

Let $\overline{x}$ be the positive-to-zero crossing point of the solid line. Then $x_0=0$ and $x_0=\overline{x}$ are the two equilibrium points of the dynamics. Since $f(x_0, x^\star(\xi), {\bf 1}) > 0$ for $0 < x_0<\overline{x}$ and $f(x_0, x^\star(\xi), {\bf 1}) < 0$ for $ \overline{x} < x_0 < m_0$, $x_0 = \overline{x}$ can easily be seen to be an asymptotically stable equilibrium of the dynamics with a basin that consists of all points but $x_0 = 0$. Thus the infection sustains so long as the initial infection level $x_0(0) > 0$.

The dashed line corresponds to the case when cross-community interactions are absent, i.e., $\nu_1=\nu_2=0$ and is the plot of the function $f(x_0, {\bf 0}, {\bf 0})$ in (\ref{eqn:zero-cross-community-upperbound}).
\end{itemize}

We now present a sufficient condition on the interaction rate parameters $\gamma_k,k=0,1,\ldots,K$ and $\nu_k,k=1,2,\ldots,K$ for the infection to sustain.

Observe that $f(x_0, x^\star(\xi), {\bf 1})=0$ when $x_0=0$. For the equilibrium $x_0=0$ to be unstable, a sufficient condition is
\begin{align}\label{eqn:cond-unstable-equilibrium}
 \left.  \frac{df}{dx_0} \right \vert_{x_0=0+}>0.
\end{align}
If the above condition is satisfied, then $f(x_0, x^\star(\xi(x_0)), {\bf 1})>0$ when $x_0$ is sufficiently close to zero. This causes $x_0$ to move away from the equilibrium point $x_0 = 0$. We now derive a sufficient condition for (\ref{eqn:cond-unstable-equilibrium}) to hold. Observe that
 \begin{align}
 \nonumber
 &\frac{df}{dx_0}=-\mu_0+ m_0 \gamma_0 - 2 \gamma_0 x_0 + \\
 \label{eqn:partial-derivative-driving-function}
 &\sum_{k=1}^{K} \left[ m_0\nu_k\frac{dx_k^{\star}(\nu_k x_0)}{dx_0}-\nu_k \left( x_0\frac{dx_k^{\star}(\nu_k x_0)}{dx_0}+x_k^{\star}(\nu_k x_0)\right) \right]
\end{align}
Differentiating (\ref{eqn:xk*}) and setting $x_0=0$, we get
\begin{align}\label{eqn:partial-derivative-equilibrium-point}
 \left.\frac{dx_1^{\star}(\nu_k x_0)}{dx_0}\right\vert_{x_0=0+}=\frac{m_k \nu_k}{b_k}.
\end{align}
Evaluating (\ref{eqn:xk*}) at $x_0=0$, we get
\begin{align}\label{eqn:equilibrium-point-evaluation}
 \left. x_k^{\star}(\nu_k x_0)\right \vert_{x_0=0}=0.
\end{align}
Substituting (\ref{eqn:partial-derivative-equilibrium-point}) and (\ref{eqn:equilibrium-point-evaluation}) in (\ref{eqn:partial-derivative-driving-function}) and evaluating (\ref{eqn:partial-derivative-driving-function}) at $x_0=0$, we get
\begin{align}
 \nonumber
 \left.  \frac{df}{dx_0} \right \vert_{x_0=0}
 & = -\mu_0 + m_0 \gamma_0 +m_0\sum_{k=1}^{K}m_k\frac{\nu_k^2}{b_k} \\
 \label{eqn:condition-sustain}
 & = -b_0 + m_0\sum_{k=1}^{K}m_k\frac{\nu_k^2}{b_k}.
\end{align}
which is strictly positive if
\[
  b_0 < m_0\sum_{k=1}^{K}m_k\frac{\nu_k^2}{b_k},
\]
which is then a sufficient condition for the infection to sustain.


\section{Optimal Control of the $(0,\infty)$-system}\label{sec:opt-control-reduced-system}

We now proceed to derive the optimal control on the $(0,\infty)$-system. Recall the dynamics of the $(0,\infty)$-system from (\ref{eqn:0-infty-equilibrium}):
\begin{align}
\dot{x}_0(t) & =  f(x_0(t), x^\star(\xi(t)), u(t)),
\end{align}
where
\begin{align*}
  f(&x_0, x^\star(\xi), u) \\
  & := - \mu_0 x_0 + (m_0 - x_0) \Big( \gamma_0 x_0 + \sum_{k=1}^K \nu_k x^{\star}_k(\xi_k) \Big) \cdot u, \\
  \xi_k & = \nu_k x_0 u, \quad k = 1, \ldots, K.
\end{align*}
Define
\begin{equation}
  \label{eqn:network-effect-R-alone}
  R(x_0 u) := \sum_{k=1}^K \nu_k x^\star_k(\nu_k x_0 u).
\end{equation}
Then we may write
\begin{align}
  \nonumber
  \dot{x}_0(t) & = - \mu_0 x_0(t) \\
  \label{eqn:network-effect}
  & ~ ~~+ (m_0 - x_0(t)) \left( \gamma_0 x_0(t) + ~R(x_0(t)u(t)) \right) \cdot u(t). \quad \quad
\end{align}
One may view $R(x_0u)$ as the \emph{effect of the network} on the evolution of the infection in the mobile population, when the infection level in the mobile population is $x_0$. Of course, when $u=0$, there is no effect of the network on the evolution of $x_0$. The dependence of $R$ on $u$ makes this problem a little more intricate. When $u=1$, by virtue of the fact that $x^\star_k(\xi_k)$ is increasing and concave in $\xi_k$, we have the following:
\begin{align}
  \nonumber
  &\mbox{The mapping $x_0 \mapsto R(x_0)$ is strictly increasing}\\
  \label{eqn:R-is-increasing-concave}
  &\mbox{and strictly concave.}
\end{align}

We now find the optimal control on the $(0,\infty)$-system.

Let us now recall the control cost on the $(0,\infty)$-system. The control variable $u(t)$ satisfies $u(t) \in [0,1]$ for all $t \in [0,T]$. Fix a deterministic measurable $u(\cdot)$ and call this open-loop policy $\pi$. The cost (\ref{eqn:cost-function}) of this open loop policy is
\begin{align}
  \nonumber
  V^{(0,\infty)}(\zeta(0), \pi) & = \int_{[0,T]} u(t) ~dt - \zeta_1(T) \\
  \label{eqn:cost}
  & = \int_{[0,T]} u(t) ~dt - x_0(T),
\end{align}
where the running cost appearing in (\ref{eqn:value-function}) is taken to be $r_1(x_0(t),u(t)) = u(t)$ and the controller experiences a terminal cost of $r_2(x_0(T)) = -x_0(T)$. In particular, more the number of infected nodes, lesser the cost, and greater the reward. As will be obvious, other scale factors for costs can be easily considered and generically they will not affect the bang-bang nature of the solution. For example, in the case of disease control, we can choose the cost function to minimize as
\[\int_{[0,T]} (1-u(t)) ~dt + x_0(T), \]
and the results can be easily adapted. Denote $J(u):=V^{(0,\infty)}(\zeta(0), \pi)$ for the sake of notational simplicity in the rest of this section.

\subsection{An Artificial Case: Network Effect Always Present}

Let us first consider a simpler problem where the network effect is always felt, even when $u=0$. By this, we mean
\begin{align}
\label{eqn:network-effect-without-u}
  \dot{x}_0(t)= - \mu_0 x_0(t) + (m_0 - x_0(t)) (\gamma_0 x_0(t) + R(x_0(t)) ) u(t).
\end{align}
The difference between (\ref{eqn:network-effect-without-u}) and (\ref{eqn:network-effect}) is that the network effect is $R(x_0)$ instead of $R(x_0 u)$. The modification has a simple solution, and provides a bound on the optimal cost for the original problem.

Observe that when the initial condition $x_0(0) = 0$, since $R(0) = 0$, the system does not rise from $0$ for any control, and it follows that $x_0(t) = 0$ for all $t \in [0,T]$. We therefore assume that $x_0(0) > 0$. In this latter case, the evolution of the system is lower bounded by the zero-control solution $x_0(0) \exp \{ - \mu_0 t \}$, i.e., for any controlled trajectory $x_0(\cdot)$, we have $x_0(t) \geq x_0(0) \exp \{ - \mu_0 t \} > 0$, for all $t \in [0,T]$.

\vspace*{.25cm}

\begin{theorem}
\label{thm:no-network-effect}
Consider the finite horizon optimal control problem for the dynamics in (\ref{eqn:network-effect-without-u}) with cost given by (\ref{eqn:cost}). Let $x_0(0) > 0$. An optimal control $u^*(\cdot)$ exists, and is a threshold policy: there exists a $\tau \in [0,T]$ such that $u^*(t) = {\bf 1} \{ t \geq \tau \}$.
\end{theorem}

\vspace*{.25cm}

\begin{IEEEproof}
We first show the existence of the optimal control for the system in (\ref{eqn:network-effect-without-u}). Consider the system in (\ref{eqn:network-effect-without-u}) with an additional variable $y(t)$ with dynamics
\begin{align}\label{eqn:augmented-variable}
\dot{y}(t)=u(t),~y(0)=0.
\end{align}
Let the cost function of this augmented system be $\tilde{J}(u)=y(T)-x_0(T)$. The system (\ref{eqn:network-effect-without-u}) and the augmented system  (\ref{eqn:network-effect-without-u})-(\ref{eqn:augmented-variable}) are equivalent in the following sense: for the same path-wise control $u(t)$ applied to  both the systems and for the same initial condition, the values of the cost functions $J(u)$ and $\tilde{J}(u)$ are equal. Hence it suffices to show the existence of optimal control for the augmented system.

Denote the reachable set for the augmented system with initial condition $x_0(0)$ as
\begin{eqnarray*}
S_T(x_0(0)) = \{(x_0(T),y(T)):x_0,y \mbox{ solves } (\ref{eqn:network-effect-without-u}) \mbox{ and } (\ref{eqn:augmented-variable}) \\
 \mbox{ for some admissible control } u(\cdot)\}.
\end{eqnarray*}
Since the cost function $\tilde{J}(u)$ depends only on the terminal state $(x_0(T),y(T))$ of the augmented system, the optimal control problem is equivalent to the problem of minimizing the cost function over the set $S_T(x_0(0))$.

We show compactness of $S_T(x_0(0))$ via Filippov's theorem \cite[p. 149-150]{liberzon2012calculus}. It is straightforward to check that solutions exist for the augmented system for every admissible control $u(\cdot)$. Denote the right-hand side of (\ref{eqn:network-effect-without-u}) as $\hat{f}(x_0,u)$. For a fixed $x_0$, the set $\{\hat{f}(x_0,u):u \in [0,1]\}$ is an affine segment in $\mathbb{R}^2$ that is convex and compact, and this is sufficient to conclude, by Filippov's theorem, that $S_T(x_0(0))$ is compact. This guarantees the existence of a minimizer of $\tilde{J}(u)$ in $S_T(x_0(0))$ and the control corresponding to the minimizer is optimal.

We will now show that the optimal control is necessarily of a threshold nature via Pontryagin's minimum principle. Consider the original non-augmented system (\ref{eqn:network-effect-without-u}) with cost given by (\ref{eqn:cost}). By defining
\begin{eqnarray}
  \label{eqn:alpha}
  \alpha(x_0) & := & - \mu_0 x_0 \\
  \label{eqn:beta}
  \beta(x_0)  & := & (m_0 - x_0) ( \gamma_0 x_0 + R(x_0)),
\end{eqnarray}
the state evolution equation (\ref{eqn:network-effect-without-u}) can we written as
\[
  \dot{x_0}(t) = \alpha(x_0(t)) + \beta(x_0(t)) \cdot u(t).
\]
Use $p$ to denote the co-state. The Hamiltonian for the system is
\begin{eqnarray}
  \nonumber
  H(x_0, u, p) & := & u + p[\alpha(x_0) + \beta(x_0) \cdot u] \\
  \nonumber
  & = & p\alpha(x_0) + [1 + p\beta(x_0)] \cdot u \\
  \label{eqn:Hamitonian-phi}
  & = & p\alpha(x_0) + \phi \cdot u
\end{eqnarray}
where $\phi := [1 + p\beta(x_0)]$ is the so-called {\em switching function}.

Pontryagin's minimum principle for the fixed time horizon, free terminal state, but with terminal cost, is the following (see for example \cite[Ch. 7, Prop. 3.3.1]{bertsekas1995dynamic}). Let $x_0^*(t), u^\star(t), p(t)$ be the optimal trajectories of the state, control, and the corresponding co-state variables, respectively. Then
\begin{itemize}
  \item The optimal state evolution is given by $$\dot{x}^\star_0(t) = \alpha(x^\star_0(t)) + \beta(x^\star_0(t)) \cdot u^\star(t),$$ with initial condition $x^\star_0(0) = x_0(0)$ given.
  \item The co-state evolution is given by
  \begin{align}
    \dot{p}(t) = & - \frac{\partial}{\partial x_0} H(x^\star_0(t), u^\star(t), p^\star(t)) \nonumber\\
    = & - p(t) \left( \frac{\partial}{\partial x_0}\alpha(x_0(t)) + \frac{\partial}{\partial x_0}\beta(x_0(t)) \cdot u^\star(t) \right),\label{eqn:costate-evolution} \\
    p(T) = & \frac{\partial}{\partial x_0}(-x_0) = -1,\nonumber
  \end{align}
  where the boundary condition is fixed by the terminal cost $-x_0(T)$.
  \item For each $t \in [0,T]$, we have $H(x_0^\star(t), u^\star(t), p(t)) \leq H(x_0^\star(t), u, p(t)) \quad \forall u \in [0,1]$.
  \item There is a constant $c$ such that $H(x_0^\star(t), u^\star(t), p(t)) = c$ for all $t \in [0,T]$.
\end{itemize}

Let us now deduce some facts about the optimal control.

(a) By the minimality criterion for the control variable, the third bullet above, and from the affine dependence of the Hamiltonian in $u$, as can be seen in (\ref{eqn:Hamitonian-phi}), we must have
\begin{equation}
  \label{eqn:u*}
  u^\star(t) = \left\{
    \begin{array}{ll}
      1 & \mbox{if } \phi(t) < 0\\
      0 & \mbox{if } \phi(t) > 0\\
      \Box & \mbox{if } \phi(t) = 0.
    \end{array}
  \right.
\end{equation}
In b) below, we will argue that there is at most one point of time where such a switch happens, and hence the value ``$\Box$'' in (\ref{eqn:u*}) is inconsequential.

b) The switching function $\phi(t)$ is strictly decreasing in time. Details follow.

If we take the time derivative of the switching function, we get
\begin{align}
  \nonumber
  \dot{\phi}(t)
  = & ~ \dot{p}(t) \beta(x_0^*(t)) + p(t) \frac{\partial}{\partial x_0} \beta(x_0^*(t)) \dot{x}^*_0(t) \\
  \nonumber
  = & - p(t) \Big( \frac{\partial}{\partial x_0}\alpha(x_0(t)) + \frac{\partial}{\partial x_0}\beta(x_0(t)) \cdot u^*(t) \Big) \beta(x_0^*(t)) \\
  \nonumber
  & + p(t) \frac{\partial}{\partial x_0} \beta(x_0^*(t))
  \left[ \alpha(x^*_0(t)) + \beta(x^*_0(t)) \cdot u^*(t) \right] \\
  \label{eqn:switching-function-evolution}
  = & ~ p(t) \cdot [\alpha, \beta](x^*_0(t)),
\end{align}
where
\[
  [\alpha, \beta](x_0) := \alpha(x_0) \frac{\partial}{\partial x_0} \beta(x_0) -  \beta(x_0) \frac{\partial}{\partial x_0} \alpha(x_0)
\]
is the Lie bracket of the two differentiable functions $\alpha$ and $\beta$ (in general two vector fields).

Next, observe from (\ref{eqn:costate-evolution}) that the co-state $p(t)$ can never take the value zero for any time in the time horizon $[0,T]$; otherwise, the terminal value of $-1$ will not be reached. Indeed, the co-state cannot even change sign; otherwise, by continuity of the co-state evolution in time, a zero value is attained at some time during $[0,T]$, and we just ruled this out. Thus $p(t) < 0$ for each $t$ in $[0,T]$.

We next claim that the Lie bracket satisfies $[\alpha, \beta](x_0) > 0$ for all feasible $x_0 \in (0, m_0]$. In particular, $[\alpha, \beta](x^*_0(t)) > 0$ for all $t$ in $[0,T]$, and so, from (\ref{eqn:switching-function-evolution}) and the fact that $p(t) < 0$ for all $t \in [0,T]$, we have that $\phi(\cdot)$ is a strictly decreasing function of time.

We now prove the claim $[\alpha, \beta](x_0) > 0$ for all $0 < x_0 \leq m_0$. Since
\[
  [\alpha, \beta](x_0) = (\alpha(x_0))^2 \cdot \frac{\partial}{\partial x_0} \left( \frac{\beta(x_0)}{\alpha(x_0)} \right),
\]
we will show that the second term is strictly positive, or equivalently, $\beta/\alpha$ is strictly \emph{increasing} in $x_0$. Let us write $\beta(x_0) = (m_0 - x_0) \tilde{R}(x_0)$, where $\tilde{R}(x_0) := \gamma_0 x_0 + R(x_0)$. Then, for $x_0 > 0$, we have
\begin{eqnarray*}
  \frac{\beta(x_0)}{\alpha(x_0)} & = & \frac{(m_0 - x_0) \tilde{R}(x_0)}{-\mu_0 x_0} \\
   & = & - \frac{1}{\mu_0} \left[ m_0 \tilde{R}(x_0) / x_0 - \tilde{R}(x_0) \right],
\end{eqnarray*}
so that it suffices to show that the term within square brackets is strictly \emph{decreasing} in $x_0$. Its derivative with respect to $x_0$ is
\begin{eqnarray*}
  \lefteqn{\left[ m_0 \tilde{R}(x_0) / x_0 - \tilde{R}(x_0) \right]'} \\
  & = &  m_0 \left( \frac{x_0 \tilde{R}'(x_0) - \tilde{R}(x_0)}{x_0^2} \right) - \tilde{R}'(x_0) \\
  & = & \frac{1}{x_0} \left( (m_0 - x_0) \tilde{R}'(x_0) - \frac{m_0\tilde{R}(x_0)}{x_0} \right) \\
  & = & \frac{m_0}{x_0} \left( \left(1 - \frac{x_0}{m_0}\right) \tilde{R}'(x_0) - \frac{\tilde{R}(x_0)}{x_0} \right) \\
  & < & \frac{m_0}{x_0} \left(\tilde{R}'(x_0) - \frac{\tilde{R}(x_0)}{x_0} \right) \\
  & < & 0.
\end{eqnarray*}
The penultimate inequality follows because $\tilde{R}(x_0) = \gamma_0 x_0 + R(x_0)$ is strictly increasing in $x_0$ for the following reasons. First, the network effect $R(x_0)$ has this property by (\ref{eqn:R-is-increasing-concave}), and hence $\tilde{R}'(x_0)$ is strictly positive. Second, dropping the factor $(1 - x_0 / m_0)$ will result in a strict increase since $x_0 > 0$. The last inequality follows because $\tilde{R}(x_0)$ is strictly concave in $x_0$, because the network effect $R(x_0)$ has this property and $\tilde{R}(x_0)$ is an affine modification of $R(x_0)$. This completes the proof of the claim that $[\alpha, \beta](x_0) > 0$ for all $x_0 > 0$.

We have thus established that the switching function is strictly decreasing.

c) Since $\phi(t)$ is strictly decreasing, $\phi(t)$ can take the value zero at not more than one point, say at $t = \tau$. It follows from (\ref{eqn:u*}) that $u^*(t)$ is increasing, and the optimal policy is a threshold policy given by ${\bf 1}\{t \geq \tau\}$ for some $\tau \in [0,T]$, where we have chosen ``$\star$'' value in (\ref{eqn:u*}) to be 1 so that the control is right-continuous with left limits. If $\tau = T$, then $u^*(t) = 0$ for all $t \in [0,T)$.

This completes the proof of Theorem \ref{thm:no-network-effect}.
\end{IEEEproof}

\subsection{Network Effect Modulated by the Control Variable}

Let us now consider the original dynamics as given in (\ref{eqn:network-effect}) when the network effect is modulated by the control variable $u$. Immediately after the proof, we highlight the reason for bringing in the artificial system with no network effect.

\emph{Proof of Theorem \ref{thm:opt-control}}: Let $x$ and $\hat{x}$ denote the solution to dynamics (\ref{eqn:network-effect}) and (\ref{eqn:network-effect-without-u} ) respectively when the same path-wise control $u(t)$ is applied to both the systems. Assume $x(0)=\hat{x}(0)$. Since $R(\cdot)$ is increasing, we have
\begin{align}\label{eqn:R-increasing}
R(x_0u) \leq R(x_0),~\forall ~ u \in [0,1],
\end{align}
Denote the right-hand side of (\ref{eqn:network-effect}) as $f(x_0,u)$ and the right-hand side of (\ref{eqn:network-effect-without-u}) as $\hat{f}(x_0,u)$. As a result of (\ref{eqn:R-increasing}), we have $f(x_0,u)\leq \hat{f}(x_0,u)$, and a quick examination of the two quantities yields that equality holds if and only if $u \in \{0,1\}$. Since $f(x_0,u)\leq \hat{f}(x_0,u)$, we have $x(t)\leq\hat{x}(t),~ \forall t \in [0,T]$. Since $x(T)\leq\hat{x}(T)$, the cost incurred for the  actual system (\ref{eqn:network-effect}) is larger than that of the artificial system in (\ref{eqn:network-effect-without-u}) when the same path-wise control is applied to both the systems. This shows that the optimal cost for the actual system in (\ref{eqn:network-effect}) is greater than or equal to the optimal cost of the artificial system in (\ref{eqn:network-effect-without-u}).

Further, it is easy to see that $f(x_0,u)=\hat{f}(x_0,u)$ for $u \in \{0,1\}$. Therefore, $x(T)=\hat{x}(T)$ when $u=u^{\star}$. So the cost incurred is equal for both systems when $u^{\star}$, the optimal control for the artificial system, is applied. Hence $u^{\star}(t)$ must also be the optimal control for the actual system (\ref{eqn:network-effect}).

This completes the proof of Theorem \ref{thm:opt-control}.$\hfill \IEEEQEDclosed$

\vspace*{.2cm}

Some remarks on our strategy for the proof of Theorem \ref{thm:opt-control} are in order.  Our strategy was to use the artificial system (\ref{eqn:network-effect-without-u}) to lower bound the cost of the actual system (\ref{eqn:network-effect}), and then show that this lower bound is attained. The reason for this indirect approach is that the set $\{ f(x_0, u) : u \in [0,1] \}$, where $f(x_0, u)$ is the right-hand side of (\ref{eqn:network-effect}), is not convex. One cannot then directly apply Filippov's theorem to conclude the existence of an optimal control. However, Filippov's theorem is indeed applicable for the artificial system since the set $\{ \hat{f}(x_0, u) : u \in [0,1] \}$, where $\hat{f}(x_0, u)$ is the right-hand side of (\ref{eqn:network-effect-without-u}), is convex and compact.

\section{Conclusion}
\label{sec:conclusion}
In this paper, we considered an information diffusion problem in a network with $K$ (almost) isolated communities and a mobile community. We studied the evolution of the infection level in the system by analyzing the evolution of the empirical measure of members possessing the information, $\zeta^{(\varepsilon,n)}(t)$. We derived the fluid limit of the original population model as the population size $n$ is scaled to infinity. The fluid limit model is a  two timescale dynamical system with a timescale separation parameter $\varepsilon$. We then obtained a dynamical system of reduced dimensions by driving the parameter $\varepsilon$ to zero. We studied the infection spread in the reduced dynamical system and showed the following.
\begin{itemize}
 \item The evolution of the infection level in the original population model and the reduced dynamical system are close to each other with high probability.
 \item The control that minimizes the cost function $\int_0^Tu(t)-\zeta_0(T)$ on the reduced dynamical system is a bang-bang control as given in (\ref{eqn:optimal-control}).
 \item The bang-bang control is nearly optimal for the original population model among those controls that have, outside a set of probability zero, control process sample paths with at most some specified $B$ number of discontinuities.
 \end{itemize}

\begin{figure}
\centering
  \includegraphics[scale=0.65]{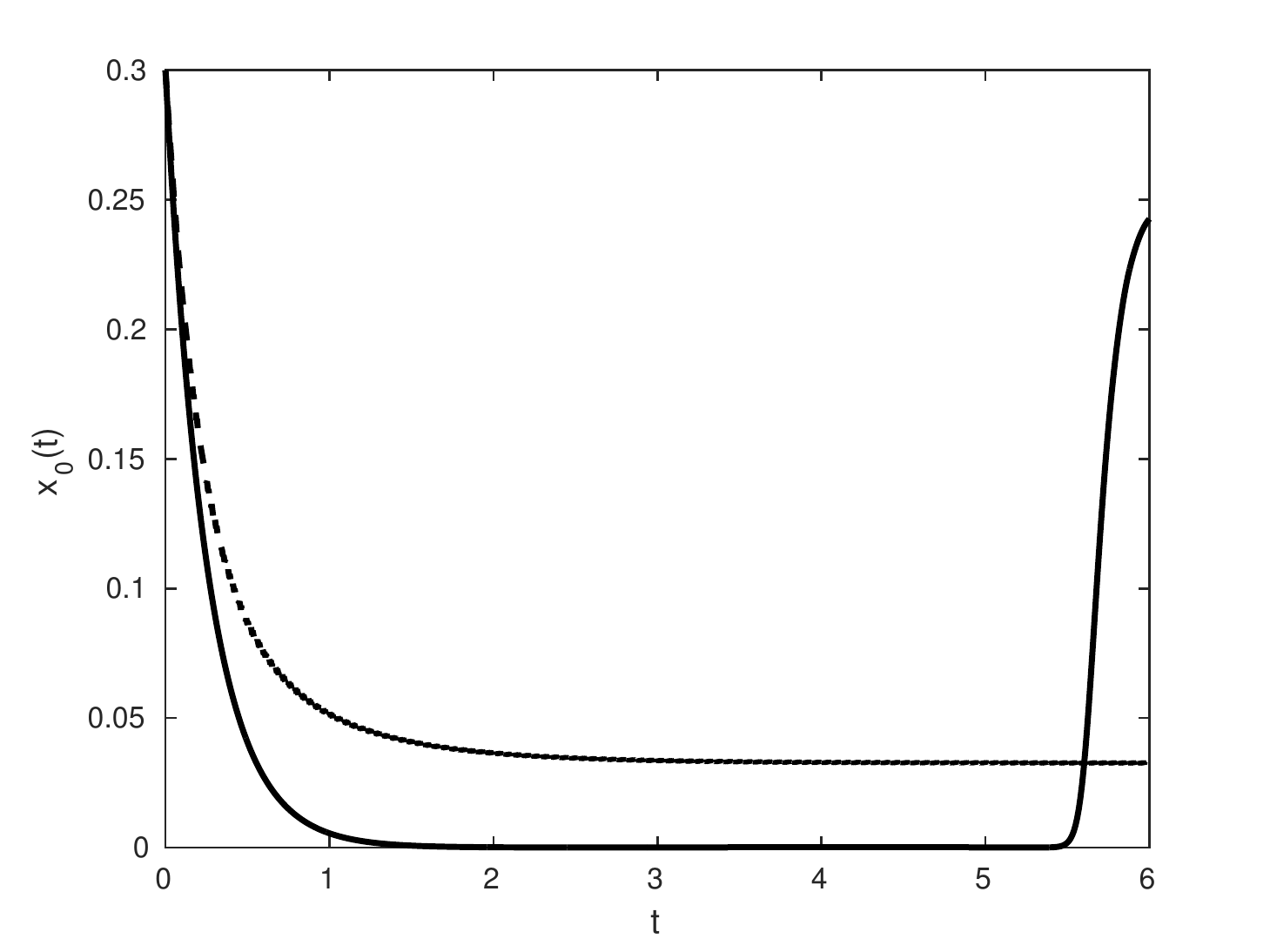}\\
  \caption{Comparison of a rapidly switched control policy with the optimal policy for a population model with two isolated communities and a mobile community. The rapidly switched control policy is a square wave that alternates between $0$ and $1$ at a rate of 100 times per second and has the same duty cycle as the optimal policy.}
  \label{fig:robust-control}
\end{figure}

Though the bang-bang control is optimal for the cost function considered in this paper, it is not robust. In the reduced dynamical system, the infection level in the mobile community $x_0(t)$ decays exponentially when the control $u^{\star}(t)$ is zero. The solid line in Figure \ref{fig:robust-control} depicts this situation. While it never dies out in the limiting reduced dynamical system, in the actual finite population model, there is the possibility that it may die out and reach the absorbing state of no infection before the threshold time $\tau$ when $u^{\star}(t)$ is switched to $1$. This is because our optimization criterion is minimization of average cost, not maximization of the probability that infection sustains.

If the criterion is maximization of the probability that infection sustains, the dashed line in Figure \ref{fig:robust-control} shows the evolution of $x_0(t)$ when an alternative relaxed control is applied. The relaxed control rapidly switches between on and off but maintains the same duty cycle as the optimal control policy over the finite time horizon $T$. This helps in keeping the infection level above a certain threshold and increases the probability of survival of the infection. The optimal control for maximization of probability of survival of infection or minimization of cost under the constraint that the infection level stays above a certain threshold can be a topic for future research.

\appendices

\section{Proof of properties of the component functions $x_k^*(\xi_k)$ and $\xi_k x^*_k(\xi_k)$}
\label{app:equilibrium-point-concavity-convexity}

\begin{prop}\label{prop:equilibrium-point-concavity}
The mapping $\xi_k \mapsto x^*_k(\xi_k)$ is a strictly increasing and strictly concave function.
\end{prop}
\begin{IEEEproof}
For ease of exposition, we do not indicate the subscript $k$. From (\ref{eqn:xk*}), it suffices to show that
\[
  \xi \mapsto r(\xi) := - (b + \xi) + \sqrt{(b+\xi)^2 + 4 \gamma m \xi}
\]
is strictly increasing and strictly concave. Taking derivative, we get
\begin{eqnarray}
  \nonumber
  r'(\xi) & = & -1 + \frac{(b+\xi) + 2 \gamma m}{\sqrt{(b+\xi)^2 + 4 \gamma m \xi}} \\
  \label{eqn:r'}
  & = & \frac{A}{B} - 1
\end{eqnarray}
where $A:= (b+\xi) + 2 \gamma m$ and $B:=\sqrt{(b+\xi)^2 + 4 \gamma m \xi}$. Now
\begin{eqnarray}
  \nonumber
  A^2 - B^2
  & = & \left( (b + \xi) + 2 \gamma m \right)^2 - \left( (b+\xi)^2 + 4 \gamma m \xi \right) \\
  \nonumber
  & = & (b + \xi)^2 + 4 \gamma^2 m^2 + 4 \gamma m (b + \xi) \\
  \nonumber
  & & -  \left( (b+\xi)^2 + 4 \gamma m \xi \right) \\
  \nonumber
  & = & 4 \gamma m (\gamma m + b) \\
  \label{eqn:A2-B2}
  & = & 4 \gamma m \mu \hspace*{1cm} (\mbox{since } b = \mu - m \gamma)\\
  \label{eqn:A2-B2ispositive}
  & > & 0.
\end{eqnarray}
Hence $A^2 > B^2$, and so $A/B > 1$. Substituting this in (\ref{eqn:r'}), we get $r'(\xi) > 0$, and this establishes the strictly increasing property.

To show strict concavity, from (\ref{eqn:r'}), we get
\begin{eqnarray}
  \nonumber
  r''(\xi) & = & \frac{B \frac{\partial A}{\partial \xi} - A \frac{\partial B}{\partial \xi}}{B^2} \\
  \nonumber
  & = & \frac{1}{B^2} \left( B - \frac{A^2}{B} \right) \\
  \label{eqn:r''}
  & = & \frac{B^2 - A^2}{B^3} \\
  \nonumber
  & < & 0,
\end{eqnarray}
where the last strict inequality follows from (\ref{eqn:A2-B2ispositive}). Hence $r(\xi)$ is strictly concave.
\end{IEEEproof}
\begin{prop}\label{prop:equilibrium-point-convexity}
The mapping $\xi_k \mapsto \xi_k x^*_k(\xi_k)$ is a strictly increasing and a strictly convex function.
\end{prop}
\begin{IEEEproof}
To show that $\xi x^*(\xi)$ is strictly increasing and strictly convex, as above, it suffices to show that $\xi r(\xi)$ is strictly increasing and strictly convex. Strictly increasing property follows immediately because both component functions $\xi$ and $r(\xi)$ that make up the product are strictly increasing.

We now show $(\xi r(\xi))'' = 2 r'(\xi) + \xi r''(\xi)$ is strictly positive. From the formulas for $r'$ and $r''$ in (\ref{eqn:r'}) and (\ref{eqn:r''}), respectively, we get
\begin{eqnarray*}
  (\xi r(\xi))'' & = & 2 \left( \frac{A}{B} - 1 \right) + \xi \left( \frac{B^2 - A^2}{B^3} \right) \\
  & = & \frac{2(A-B)}{B^3} \left[ B^2 - \xi \left( \frac{A+B}{2} \right) \right] \\
  & > & 0.
\end{eqnarray*}
The last inequality follows from $A > B$ and the stronger inequality
\[
  B^2 > \xi A > \xi \left( \frac{A + B}{2} \right).
\]
The second of these inequalities holds simply because $A > B$. The first of these is obtained as follows:
\begin{eqnarray*}
  B^2 - \xi A
  & = & (b + \xi)^2 + 4 \gamma m \xi - \xi [(b+\xi) + 2 \gamma m] \\
  & = & (b+\xi) b + 2 \gamma m \xi > 0.
\end{eqnarray*}
This completes the proof of the strict convexity of $\xi r(\xi)$, and hence of the strict convexity of $\xi x^*(\xi)$.
\end{IEEEproof}

\section{A general mean-field convergence result for a closed and controlled system}\label{app:appendix-general-convergence}
In this section, we state and prove a more general and refined version of Theorem \ref{thm:cost-equality-n-infty}.
The following assertions hold true for the transition rate matrices $\Lambda^{(n)}$ of the $(\varepsilon,n)$-system and $\Lambda$ of the $(\varepsilon,\infty)$-system.

Below, the quantities $\zeta,\zeta' \in \mathcal{E}$  and $u \in [0,1]$.
\begin{itemize}
  \item[({\bf F1})] The functions $\lambda_{i,j}^n \rightarrow \lambda_{i,j}$ uniformly as $n \rightarrow \infty$. More precisely,
  \[
    \sup_{\zeta, u} |\lambda_{i,j}^n(\zeta,u) - \lambda_{i,j}(\zeta,u)| \leq I(n) \ra 0 \mbox{ as } n \ra \infty.
  \]

  \item[({\bf F2})] The mappings $\zeta \mapsto \lambda_{i,j}(\zeta, u)$ are uniformly Lipschitz continuous, uniformly in $u \in [0,1]$. More precisely, there is an $L_2$ such that for every $\zeta, \zeta', u$, we have
  \[
    |\lambda_{i,j}(\zeta', u) - \lambda_{i,j}(\zeta, u)| \leq L_2 ||\zeta - \zeta'||.
  \]

\end{itemize}

The following bounds will be needed. Assertion ({\bf F1}) implies that there is a suitable constant $L_3$ such that, for all $\zeta$ and $u$, we have
\begin{equation}
  \label{eqn:A1-consequence-1}
  ||\left[\Lambda^{(n)}(\zeta,u) - \Lambda(\zeta,u))\right]^* \zeta|| \leq L_3 I(n).
\end{equation}
Moreover, the constant can be chosen so that, for all $n$, $\zeta$, $u$, we have
\begin{equation}
  \label{eqn:A1-consequence-2}
  ||\left[\Lambda^{(n)}(\xi,u)\right]^* \zeta|| \leq L_3.
\end{equation}
Assertion ({\bf F2}) implies that there is a suitable constant $L_2'$ such that, for all $\zeta, \zeta', u$, we have
\begin{eqnarray}
  \label{eqn:A2-consequence-1}
  ||\left[\Lambda(\zeta,u)^* (\zeta - \zeta')\right]|| & \leq & L_2' ||\zeta - \zeta'||, \\
  \label{eqn:A2-consequence-2}
  ||\left[\Lambda(\zeta,u) - \Lambda(\zeta',u)\right]^* \zeta || & \leq & L_2' ||\zeta - \zeta'||.
\end{eqnarray}

Following the idea of \cite{gast2012mean}, we couple the dynamics of the $(\varepsilon,n)$-system with that of the limiting system by employing the same path-wise control. Let $\pi$ denote any non-anticipative control policy on the $(\varepsilon,n)$-system. It could be a feedback policy on $(\varepsilon,n)$-system or a simple open-loop policy. Let $\zeta^{(\varepsilon,n)}(t)$ be the empirical distribution of the $(\varepsilon,n)$-system at time $t$ and let $u(t)$ be the resulting control action; $\zeta^{(\varepsilon,\infty)}(\cdot)$ is then the solution to (\ref{eqn:fluidlimit}) when driven by the sample path $u(\cdot)$. The evolution of the limiting system may be random on account of the possible randomness in the control. We simplify notation by ignoring $\varepsilon$ to denote $\zeta^{(\varepsilon,n)}(t)$ as $\zeta^n(t)$ and $\zeta^{(\varepsilon,\infty)}(t)$ as $\zeta^{\infty}(t)$ in the rest of the section. Let $\vert \vert \zeta^n -\zeta^\infty\vert \vert_T=\sup_{t \in [0,T]} \vert \vert \zeta^n(t) -\zeta^\infty (t)\vert \vert$.

\vspace*{.2cm}

\begin{theorem}\label{thm:mean-field-limit}
Fix $n$, a policy $\pi$ on the $(\varepsilon,n)$-system, and $b > 0$. Fix initial conditions $\zeta^{n}(0)$ and $\zeta^{\infty}(0)$ on the $(\varepsilon,n)$-system and on the limiting $(\varepsilon,\infty)$-system, respectively. Let $u(\cdot)$ be the control process sample path on the $(\varepsilon, n)$-system. Apply this control process sample path on the $(\varepsilon, \infty)$-system and call the resulting policy also as $\pi$ but on the $(\varepsilon, \infty)$-system. We have
\begin{eqnarray}
  \nonumber
  \lefteqn{ P^{\pi} \Big\{ ||\zeta^{n} - \zeta^{\infty}||_T >
  \left( ||\zeta^{n}(0) -  \zeta^{\infty}(0)|| + L_3 I(n) T + b \right)  } \\
  \nonumber
  & & \hspace*{5.5cm} \cdot \exp\{ 2L_2' T \} \Big\} \\
  & & \leq \frac{4|\mS|^2 L_3 T}{b^2 n}.
  \label{eqn:p-bound}
\end{eqnarray}
Moreover,
\begin{eqnarray}
  \lefteqn{ \left|
  V^{(\varepsilon,n)}(\zeta^n(0), \pi) - V^{(\varepsilon,\infty)}(\zeta^{\infty}(0), \pi) \right| } \nonumber \\
  \nonumber
  & \leq & L_1(T+1) \Big[  \left(|| \zeta^{n}(0) -  \zeta^{\infty}(0) || + L_3 I(n) T + b \right)\\
  \nonumber
  & & \hspace*{5.3cm} \cdot \exp\{2L_2' T\}  \\
  & & \hspace*{1.6cm} + \frac{8 |\mS|^2 L_3 T}{b^2 n} \Big]
  \label{eqn:v-bound}
\end{eqnarray}
\end{theorem}

\begin{IEEEproof}
We will use Gronwall's inequality (see e.g., \cite[Appendix B]{borkar2009stochastic}) as is customary in such proofs. Using (\ref{eqn:empdisofnthsubsystem}) and (\ref{eqn:fluidlimit}), we get
\begin{eqnarray}
  \zeta^n(t) - \zeta^{\infty}(t) & = & \zeta^{n}(0) - \zeta^{\infty}(0) \nonumber \\
  & & + ~ \int_0^t \left[ \left[\Lambda^{(n)}(\zeta^{n}(s), u(s))\right]^* \zeta^{n}(s) \right] ds \nonumber \\
  & & - \int_0^t \left[\left[ \Lambda(\zeta^{\infty}(s), u(s))\right]^* \zeta^{\infty}(s) \right] ds \nonumber \\
  \label{eqn:Mt-difference-expansion}
  & & + ~ \mathcal{M}^{n}(t),
\end{eqnarray}
where $\mathcal{M}^n(t)$ is a vector-valued square-integrable $(\F_t)$-measurable martingale.
The second term on the right-hand side of (\ref{eqn:Mt-difference-expansion}) can be manipulated, by adding and subtracting terms and by using the triangle inequality, to get
\begin{eqnarray*}
  \lefteqn{ \left\| \left[\Lambda^{(n)}(\zeta^{n}(s), u(s))\right]^* \zeta^{n}(s) - \left[\Lambda(\zeta^{\infty}(s), u(s))\right]^* \zeta^{\infty}(s)\right\| } \\
  & \leq & \hspace*{-.15cm} \left\| \left[\Lambda^{(n)}(\zeta^{n}(s), u(s))\right]^* \zeta^n(s) - \left[\Lambda(\zeta^n(s), u(s))\right]^* \zeta^n(s)\right\| \\
  & & \hspace*{-.15cm} + \|\left[ \Lambda(\zeta^n(s), u(s))\right]^* \zeta^n(s) - \left[\Lambda(\zeta^n(s), u(s))\right]^* \zeta^{\infty}(s) \| \\
  & & \hspace*{-.15cm} + \|\left[ \Lambda(\zeta^n(s), u(s))\right]^* \zeta^{\infty}(s) - \left[\Lambda(\zeta^{\infty}(s), u(s))\right]^* \zeta^{\infty}(s) \| \\
  & \leq & \hspace*{-.15cm}  L_3 I(n) + L_2'||\zeta^n(s) - \zeta^{\infty}(s)|| + L_2'||\zeta^n(s) - \zeta^{\infty}(s)||,
\end{eqnarray*}
where, to get the last inequality, we used (\ref{eqn:A1-consequence-1}), (\ref{eqn:A2-consequence-1}), and (\ref{eqn:A2-consequence-2}), which are consequences of Assumptions ({\bf F1}) and ({\bf F2}). Observe that we crucially use the fact that the employed control path is the \emph{same} in both systems. Substituting the above inequality back in (\ref{eqn:Mt-difference-expansion}), we get
\begin{eqnarray*}
  \lefteqn{\| \zeta^n(t) - \zeta^{\infty}(t) \|} \\
   & \leq & \hspace*{-.2cm}\| \zeta^n(0) - \zeta^\infty (0) \| \\
   & & \hspace*{-.2cm} + \int_0^t [L_3 I(n) + 2 L_2' \| \zeta^n(s) - \zeta^\infty(s) \|]~ds + \| \mathcal{M}^n(t) \|.
\end{eqnarray*}
Under the event $ G = \{ \| \mathcal{M}^n(\cdot) \|_T \leq b \}$, by Gronwall's lemma \cite[Appendix B]{borkar2009stochastic}, we get
\begin{eqnarray*}
  \lefteqn{ \| \zeta^n(t) - \zeta^\infty(t) \| } \\
  & \leq & (\| \zeta^n(0) - \zeta^\infty(0) \| + L_3 I(n) t + b ) \exp\{2L_2't\},
\end{eqnarray*}
and furthermore,
\begin{equation}
  \label{eqn:G-event-upper-bound}
  \| \zeta^n - \zeta^\infty \|_T \leq (\| \zeta^n(0) - \zeta^\infty(0) \| + L_3 I(n) T + b ) \exp\{2L_2'T\},
\end{equation}
so that
\begin{eqnarray}
  \lefteqn{ P^{\pi} \left\{ \| \zeta^n - \zeta^\infty \|_T > (\| \zeta^n(0) - \zeta^\infty(0) \| + L_3 I(n) T + b ) \right. } \nonumber \\
  & & \left. \hspace*{5.2cm}\cdot \exp\{2L_2'T\} \right\}  \nonumber \\
  & \leq & P^{\pi} \{ G^c \}  \nonumber \\
  & = & P^{\pi} \{ \| \mathcal{M}^n(\cdot) \|_T > b \} \nonumber \\
  & \leq & \sum_{i \in \mS} \frac{4 \mathbb{E}^{\pi} [  |\mathcal{M}^n_i(T)|^2 ] \cdot |\mS|}{b^2} \label{eqn:doob} \\
  & = &  \frac{4|\mS|^2}{b^2 n} \max_i \int_0^T \hspace*{-.2cm} \Big( \Big[\Lambda^{(n)}(\zeta^n(s), u(s))\Big]^* \zeta^n(s) \Big)_i ds ~~~\label{eqn:doob-2} \\
  \label{eqn:prob-G}
  & \leq & \frac{4|\mS|^2 L_3 T}{b^2 n},
\end{eqnarray}
where inequality (\ref{eqn:doob}) follows from Markov's inequality, the union bound, and Doob's inequality \cite[Cor.~2.17, p.~64]{ethier2009markov}. Inequality (\ref{eqn:doob-2}) follows because $n\mathcal{M}_i^n(\cdot)$ is a sum of time-inhomogeneous Poisson point processes. The variance of $n\mathcal{M}_i^n(T)$ is $n$ times the integral of the intensity over the duration of the process. Inequality (\ref{eqn:prob-G}) follows from (\ref{eqn:A1-consequence-2}). This establishes (\ref{eqn:p-bound}).

Let us now turn to (\ref{eqn:v-bound}). Observe that
\begin{eqnarray*}
  \lefteqn{ \left|V^{(\varepsilon,n)}(\zeta^n(0), \pi) - V^{\varepsilon,\infty}(\zeta^\infty(0), \pi) \right| } \\
  & \leq & \mathbb{E}^{\pi} \Big[ \int_0^T [ r_1(\zeta^n(s), u(s)) - r_1(\zeta^{\infty}(s), u(s)) ] ~ds \\
  & & \hspace*{3cm} + r_2(\zeta^n(T)) - r_2(\zeta^\infty(T)) \Big] \\
  & \leq & L_1 \mathbb{E}^{\pi} \Big[ \int_0^T \hspace*{-.15cm}\| \zeta^n(s) - \zeta^\infty(s) \| ~ds + \| \zeta^n(T) - \zeta^\infty(T) \| \Big] \\
  & \leq & L_1 (T + 1) \mathbb{E}^{\pi} \left[ \| \zeta^n - \zeta^\infty \|_T  \right]
\end{eqnarray*}
The argument under the expectation above is upper bounded by $2$, and under event $G$ is upper bounded by (\ref{eqn:G-event-upper-bound}). Since $P \{ G^c \} $ itself is upper bounded by (\ref{eqn:prob-G}), the result follows from
\begin{align*}
  \lefteqn{\mathbb{E}[\| \zeta^n - \zeta^\infty \|_T] } \\
  & \leq \mathbb{E}[\| \zeta^n - \zeta^\infty \|_T {\bf 1}_G] + 2 P\{ G^c \} \\
  & \leq (\| \zeta^n(0) - \zeta^\infty(0) \| + L_3 I(n) T + b ) \exp\{2L_2'T\}  \\
  & \hspace*{5.4cm} +  \frac{8|\mS|^2 L_3 T}{b^2 n}.
\end{align*}
This concludes the proof of Theorem \ref{thm:mean-field-limit}.
\end{IEEEproof}

\vspace*{.2cm}

We now specialize the above result to prove Theorem \ref{thm:cost-equality-n-infty}.

\vspace*{.2cm}

\begin{IEEEproof}[Proof of Theorem \ref{thm:cost-equality-n-infty}]
Using (\ref{eqn:asymptotic-rates}) in (\ref{eqn:ratetranmatrix1}) and (\ref{eqn:ratetranmatrix2}), we have $\Lambda^{(n)}=\Lambda,~\forall n$. Hence $I(n)=0$. Also, by assumption, $\zeta^n(0)=\zeta^{\infty}(0)=\zeta(0)$. Hence $|| \zeta^{n}(0) -  \zeta^{\infty}(0) || = 0$. Furthermore, it is straightforward to show that $L_2,L_2^{\prime}$, and $L_3$ are $O(1/\varepsilon)$. Thus $2L_2'T=k_2/\varepsilon$ for some positive constant $k_2$ and $4|\mS|^2 L_3 T=k_3/\varepsilon$ for some positive constant $k_3$. We thus have
\begin{eqnarray}
  \nonumber
  \Big( || \zeta^{n}(0) -  \zeta^{\infty}(0) || + L_3 I(n) T + b \Big) \exp\{ 2L_2' T \} \\
  \label{eqn:orderlogn-1}
   =  b \exp\{ k_2/\varepsilon \},
\end{eqnarray}
and by setting this to be $c/\log n$, we get $b=\frac{c \exp\{ -k_2 / \varepsilon \}}{\log n }$. Plugging this into the upper bound in \eqref{eqn:p-bound} and using $4|\mS|^2 L_3 T=k_3/\varepsilon$, we get
\begin{eqnarray*}
P^{\pi} \left\{ ||\zeta^{n} - \zeta^{\infty}||_T > \frac{c}{\log n} \right\} & \leq & \frac{4|\mS|^2 L_3 T}{b^2 n} \\
& = & \frac{\frac{k_3}{\varepsilon} \exp\{2k_2/\varepsilon \}}{c^2 n/(\log n)^2} \\
& \leq & \frac{k_3 (\log n)^3}{C c^2 n^{1-2k_2/C}} \\
& = & O(1/\log n)
\end{eqnarray*}
if $C/\log n \leq \varepsilon \rightarrow 0$ with $C > 2k_2$. Thus we obtain (\ref{eqn:theorem-p-bound}). Further, plugging these choices of $b$ and $\varepsilon$ into \eqref{eqn:v-bound} we get \eqref{eqn:theorem-v-bound}. The factor $L_1(T+1)$ is just a constant since $T$ is fixed and $L_1$ is the Lipschitz constant for the cost functions $r_1$ and $r_2$ and is therefore independent of $\varepsilon$ and $n$. The constant $\bar{C}$ may be suitably chosen so that both (\ref{eqn:theorem-p-bound}) and \eqref{eqn:theorem-v-bound} hold. This completes the proof of Theorem \ref{thm:cost-equality-n-infty}.
\end{IEEEproof}

%
%
%
%
%

\bibliographystyle{IEEEtran}
\bibliography{IEEEabrv,wisl}

\begin{thebibliography}{10}
\providecommand{\url}[1]{#1}
\csname url@samestyle\endcsname
\providecommand{\newblock}{\relax}
\providecommand{\bibinfo}[2]{#2}
\providecommand{\BIBentrySTDinterwordspacing}{\spaceskip=0pt\relax}
\providecommand{\BIBentryALTinterwordstretchfactor}{4}
\providecommand{\BIBentryALTinterwordspacing}{\spaceskip=\fontdimen2\font plus
\BIBentryALTinterwordstretchfactor\fontdimen3\font minus
  \fontdimen4\font\relax}
\providecommand{\BIBforeignlanguage}[2]{{%
\expandafter\ifx\csname l@#1\endcsname\relax
\typeout{** WARNING: IEEEtran.bst: No hyphenation pattern has been}%
\typeout{** loaded for the language `#1'. Using the pattern for}%
\typeout{** the default language instead.}%
\else
\language=\csname l@#1\endcsname
\fi
#2}}
\providecommand{\BIBdecl}{\relax}
\BIBdecl

\bibitem{athreya2018simultaneous}
S.~R. Athreya, V.~S. Borkar, K.~S. Kumar, and R.~Sundaresan, ``Simultaneous
  small noise limit for singularly perturbed slow-fast coupled diffusions,''
  \emph{arXiv preprint arXiv:1810.03585}, 2018.

\bibitem{ottaviano2017optimal}
\BIBentryALTinterwordspacing
S.~Ottaviano, F.~De~Pellegrini, S.~Bonaccorsi, and P.~Van~Mieghem, ``Optimal
  curing policy for epidemic spreading over a community network with
  heterogeneous population,'' \emph{Journal of Complex Networks}, 2017.
  [Online]. Available: \url{http://dx.doi.org/10.1093/comnet/cnx060}
\BIBentrySTDinterwordspacing

\bibitem{morton1974optimal}
R.~Morton and K.~H. Wickwire, ``On the optimal control of a deterministic
  epidemic,'' \emph{Advances in Applied Probability}, vol.~6, no.~4, pp.
  622--635, 1974.

\bibitem{behncke2000optimal}
H.~Behncke, ``Optimal control of deterministic epidemics,'' \emph{Optimal
  control applications and methods}, vol.~21, no.~6, pp. 269--285, 2000.

\bibitem{karnik2012optimal}
A.~Karnik and P.~Dayama, ``Optimal control of information epidemics,'' in
  \emph{Communication Systems and Networks (COMSNETS), 2012 Fourth
  International Conference on}.\hskip 1em plus 0.5em minus 0.4em\relax IEEE,
  2012, pp. 1--7.

\bibitem{kandhway2014run}
K.~Kandhway and J.~Kuri, ``How to run a campaign: Optimal control of {SIS} and
  {SIR} information epidemics,'' \emph{Applied Mathematics and Computation},
  vol. 231, pp. 79--92, 2014.

\bibitem{altman2010optimal}
E.~Altman, A.~P. Azad, T.~Ba\c{s}ar, and F.~De~Pellegrini, ``Optimal activation
  and transmission control in delay tolerant networks,'' in \emph{INFOCOM, 2010
  Proceedings IEEE}.\hskip 1em plus 0.5em minus 0.4em\relax IEEE, 2010, pp.
  1--5.

\bibitem{khouzani2012optimal}
M.~Khouzani, S.~Sarkar, and E.~Altman, ``Optimal dissemination of security
  patches in mobile wireless networks,'' \emph{IEEE Transactions on Information
  Theory}, vol.~58, no.~7, pp. 4714--4732, 2012.

\bibitem{kumar2018influenzing}
B.~Kumar, N.~Sahasrabudhe, and S.~Moharir, ``On influencing opinion dynamics
  over finite time horizons,'' in \emph{Mathematical Theory of Networks and
  Systems}, 2018.

\bibitem{colizza2006role}
V.~Colizza, A.~Barrat, M.~Barth{\'e}lemy, and A.~Vespignani, ``The role of the
  airline transportation network in the prediction and predictability of global
  epidemics,'' \emph{Proceedings of the National Academy of Sciences of the
  United States of America}, vol. 103, no.~7, pp. 2015--2020, 2006.

\bibitem{becker1995effect}
N.~G. Becker and K.~Dietz, ``The effect of household distribution on
  transmission and control of highly infectious diseases,'' \emph{Mathematical
  biosciences}, vol. 127, no.~2, pp. 207--219, 1995.

\bibitem{ball2004stochastic}
F.~G. Ball, T.~Britton, and O.~D. Lyne, ``Stochastic multitype epidemics in a
  community of households: Estimation of threshold parameter {R} and secure
  vaccination coverage,'' \emph{Biometrika}, vol.~91, no.~2, pp. 345--362,
  2004.

\bibitem{gopalan2011random}
A.~Gopalan, S.~Banerjee, A.~K. Das, and S.~Shakkottai, ``Random mobility and
  the spread of infection,'' in \emph{INFOCOM, 2011 Proceedings IEEE}.\hskip
  1em plus 0.5em minus 0.4em\relax IEEE, 2011, pp. 999--1007.

\bibitem{de2010optimal}
F.~De~Pellegrini, E.~Altman, and T.~Ba\c{s}ar, ``Optimal monotone forwarding
  policies in delay tolerant mobile ad hoc networks with multiple classes of
  nodes,'' in \emph{Modeling and Optimization in Mobile, Ad Hoc and Wireless
  Networks (WiOpt), 2010 Proceedings of the 8th International Symposium
  on}.\hskip 1em plus 0.5em minus 0.4em\relax IEEE, 2010, pp. 497--504.

\bibitem{kandhway2016campaigning}
K.~Kandhway and J.~Kuri, ``Campaigning in heterogeneous social networks:
  Optimal control of {SI} information epidemics,'' \emph{IEEE/ACM Transactions
  on Networking}, vol.~24, no.~1, pp. 383--396, 2016.

\bibitem{kokotovic1984applications}
P.~V. Kokotovi{\'c}, ``Applications of singular perturbation techniques to
  control problems,'' \emph{SIAM review}, vol.~26, no.~4, pp. 501--550, 1984.

\bibitem{levin1954singular}
J.~Levin and N.~Levinson, ``Singular perturbations of non-linear systems of
  differential equations and an associated boundary layer equation,''
  \emph{Journal of Rational Mechanics and Analysis}, vol.~3, pp. 247--270,
  1954.

\bibitem{liberzon2012calculus}
D.~Liberzon, \emph{Calculus of Variations and Optimal Control Theory: A Concise
  Introduction}.\hskip 1em plus 0.5em minus 0.4em\relax Princeton University
  Press, 2012.

\bibitem{bertsekas1995dynamic}
D.~P. Bertsekas, \emph{Dynamic programming and optimal control}, 3rd~ed.\hskip
  1em plus 0.5em minus 0.4em\relax Belmont, MA: Athena Scientific, 2005,
  vol.~1.

\bibitem{gast2012mean}
N.~Gast, B.~Gaujal, and J.-Y. Le~Boudec, ``Mean field for markov decision
  processes: from discrete to continuous optimization,'' \emph{IEEE
  Transactions on Automatic Control}, vol.~57, no.~9, pp. 2266--2280, 2012.

\bibitem{borkar2009stochastic}
V.~S. Borkar, \emph{Stochastic approximation: {A} dynamical systems
  viewpoint}.\hskip 1em plus 0.5em minus 0.4em\relax Springer, 2009, vol.~48.

\bibitem{ethier2009markov}
S.~N. Ethier and T.~G. Kurtz, \emph{Markov processes: Characterization and
  convergence}.\hskip 1em plus 0.5em minus 0.4em\relax John Wiley \& Sons,
  2009, vol. 282.

\end{thebibliography}

\end{document}